\documentclass[12pt,a4paper]{amsart}
\theoremstyle{plain}

\usepackage[colorlinks]{hyperref}
\usepackage{enumerate, amssymb}

\advance\hoffset-5mm \advance\textwidth40mm
\input{diagrams.tex}
\diagramstyle[scriptlabels,height=8mm,width=8mm]

\def\bdi{\begin{diagram}}
\def\edi{\end{diagram}}

\theoremstyle{plain}

\newtheorem{thm}{Theorem}[section]
\newtheorem{cor}[thm]{Corollary}
\newtheorem{lem}[thm]{Lemma}
\newtheorem{prop}[thm]{Proposition}
\theoremstyle{definition}
\newtheorem{defi}[thm]{Definition}
\newtheorem{defis}[thm]{Definitions}
\newtheorem{conj}[thm]{Conjecture}
\newtheorem{conv}[thm]{Convention}
\newtheorem{nota}[thm]{Notation}
\newtheorem{rem}[thm]{Remark}
\newtheorem{rems}[thm]{Remarks}
\newtheorem{exa}[thm]{Example}
\newtheorem{exas}[thm]{Examples}
\newtheorem{prob}[thm]{Problem}
\newtheorem{probs}[thm]{Problems}
\newtheorem{ques}[thm]{Question}
\newtheorem{sit}[thm]{}

\newcommand{\Span}{ \operatorname{{\rm Span}}}
\newcommand{\Spec}{ \operatorname{{\rm Spec}}}

\newcommand{\Sing}{ \operatorname{{\rm Sing}}}

\newcommand{\Aut}{ \operatorname{{\rm Aut}}}

\newcommand{\LND}{ \operatorname{{\rm LND}}}
\newcommand{\GL}{ \operatorname{{\rm GL}}}
\newcommand{\SL}{ \operatorname{{\bf SL}}}

\newcommand{\Proj}{ \operatorname{{\rm Proj}}}

\def\reg{{\mathop{\rm reg}}}

\def\codim{\mathop{\rm codim}}

\def\lto{\longrightarrow}
\def\hto{\hookrightarrow}
\newcommand{\emb}{\hookrightarrow}

\renewcommand{\epsilon}{\varepsilon}

\def\ol{\overline}

\def\and{\quad\mbox{and}\quad}

\newcommand{\C}{\ensuremath{\mathbb{C}}}

\newcommand{\N}{\ensuremath{\mathbb{N}}}
\newcommand{\G}{\ensuremath{\mathbb{G}}}

\newcommand{\kk}[1]{\bk^{[#1]}}

\newcommand{\bX}{{\bar X}}

\newcommand{\bY}{{\bar Y}}

\newcommand{\bQ}{{\bar Q}}

\def\fa{{\mathfrak a}}
\def\fb{{\mathfrak b}}
\def\fc{{\mathfrak c}}

\newcommand{\cF}{{\ensuremath{\mathcal{F}}}}

\newcommand{\cE}{{\ensuremath{\mathcal{E}}}}

\newcommand{\cO}{{\ensuremath{\mathcal{O}}}}
\newcommand{\cC}{{\ensuremath{\mathcal{C}}}}
\newcommand{\cD}{{\ensuremath{\mathcal{D}}}}

\newcommand{\cI}{{\ensuremath{\mathcal{I}}}}

\newcommand{\cN}{{\ensuremath{\mathcal{N}}}}

\newcommand{\cX}{{\ensuremath{\mathcal{X}}}}

\newcommand{\cV}{{\ensuremath{\mathcal{V}}}}
\newcommand{\cU}{{\ensuremath{\mathcal{U}}}}

\newcommand{\p}{\partial}
\newcommand{\de}{\delta}
\newcommand{\id}{{\rm id}}

\renewcommand{\rho}{\varrho}
\newcommand{\vp}{\varphi}
\newcommand{\brho}{\bar\varrho}

\def\bals#1\eals{\begin{align*}#1\end{align*}}
\def\bal#1\eal{\begin{align}#1\end{align}}

\def\SAut{\mathop{\rm SAut}}

\def\kk{{\Bbbk}}
\def\AA{{\mathbb A}}

\def\NN{{\mathbb N}}
\def\ZZ{{\mathbb Z}}

\def\kk{{\Bbbk}}
\def\CC{{\mathbb C}}

\def\PP{{\mathbb P}}

\renewcommand{\phi}{\varphi}

\newcommand{\bnum}{\begin{enumerate}}
\newcommand{\enum}{\end{enumerate}}
\renewcommand{\emptyset}{\varnothing}

\newcommand{\lan}{\langle}
\newcommand{\ran}{\rangle}
\newcommand{\lla}{\lan\lan}
\newcommand{\rra}{\ran\ran}

\newcommand{\la}{\label}

\newcommand{\brem}{\begin{rem}}
\newcommand{\brems}{\begin{rems}}
\newcommand{\erem}{\end{rem}}
\newcommand{\erems}{\end{rems}}
\newcommand{\bprob}{\begin{prob}}
\newcommand{\eprob}{\end{prob}}
\newcommand{\bprobs}{\begin{probs}}
\newcommand{\eprobs}{\end{probs}}
\newcommand{\bques}{\begin{ques}}
\newcommand{\eques}{\end{ques}}
\newcommand{\bexa}{\begin{exa}}
\newcommand{\bexas}{\begin{exas}}
\newcommand{\eexa}{\end{exa}}
\newcommand{\eexas}{\end{exas}}
\newcommand{\bdefi}{\begin{defi}}
\newcommand{\edefi}{\end{defi}}
\newcommand{\bdefis}{\begin{defis}}
\newcommand{\edefis}{\end{defis}}
\newcommand{\bcor}{\begin{cor}}
\newcommand{\ecor}{\end{cor}}
\newcommand{\blem}{\begin{lem}}
\newcommand{\elem}{\end{lem}}
\newcommand{\bconv}{\begin{conv}}
\newcommand{\econv}{\end{conv}}
\newcommand{\bconj}{\begin{conj}}
\newcommand{\econj}{\end{conj}}
\newcommand{\bprop}{\begin{prop}}
\newcommand{\eprop}{\end{prop}}
\newcommand{\bthm}{\begin{thm}}
\newcommand{\ethm}{\end{thm}}
\newcommand{\bnota}{\begin{nota}}
\newcommand{\enota}{\end{nota}}
\newcommand{\bsit}{\begin{sit}}
\newcommand{\esit}{\end{sit}}
\newcommand{\be}{\begin{equation}}
\newcommand{\ee}{\end{equation}}
\newcommand{\bproof}{\begin{proof}}
\newcommand{\eproof}{\end{proof}}
\def\ba{\begin{array}}
\def\ea{\end{array}}

\newcommand{\nlin}{\unitlength1mm\begin{picture}(0,9.25)
                      \put(0,0.75){\line(0,1){8.5}}
                     \end{picture}}

\newcommand{\vlin}[1]{\hspace{0.75mm}\unitlength1mm\begin{picture}
(#1,0)
                      \put(0,0){\line(1,0){#1}}
                     \end{picture}\hspace{0.75mm}\rule[-3mm]{0mm}
                     {4mm}}

\def\llin{\vlin{11.5}}
\newcommand{\lin}{\vlin{8.5}}

\newcommand{\cou}[2]{\unitlength1mm\begin{picture}(0,8)
   \put(0,0){\circle{1.5}}
   \put(0,3){\makebox(0,5)[b]{$#1$}}
   \put(0,-7){\makebox(0,4)[t]{$#2$}}
     \end{picture}
     \rule[-7mm]{0mm}{7mm}}

\newcommand{\crl}[2]{\unitlength1mm\begin{picture}(0,8)
   \put(0,0){\circle{1.5}}
   \put(-5,0){\makebox(0,5)[b]{$#1$}}
  \put(5,0){\makebox(0,5)[b]{$#2$}}
     \end{picture}
     \rule[-7mm]{0mm}{7mm}}

\newcommand{\cshiftup}[2]{\unitlength1mm\begin{picture}(0,9.25)
                      \put(0,10){\crl{#1}{#2}}
                     \end{picture}}

\thanks{The  research of the second author
was supported by an NSA grant No. \ H98230-10-1-0185. This work
was done during a stay of the authors at the Max Planck Institut
f\"ur Mathematik at Bonn, and a stay of the first and the second
authors at the Institut Fourier, Grenoble. The authors thank these
institutions for hospitality.}

\begin{document}
\title[Gromov-Winkelmann Theorem]{The Gromov-Winkelmann Theorem for
flexible varieties}

\author{H.\ Flenner, S.\ Kaliman,
M.\ Zaidenberg}
\address{Fakult\"at f\"ur
Mathematik, Ruhr Universit\"at Bochum, Geb.\ NA 2/72,
Universit\"ats\-str.\ 150, 44780 Bochum, Germany}
\email{Hubert.Flenner@rub.de}
\address{Department of Mathematics,
University of Miami, Coral Gables, FL 33124, USA}
\email{kaliman@math.miami.edu}
\address{Universit\'e Grenoble I, Institut Fourier, UMR 5582
CNRS-UJF, BP 74, 38402 Saint Martin d'H\`eres c\'edex, France}
\email{Mikhail.Zaidenberg@ujf-grenoble.fr}

\begin{abstract} An affine variety $X$
of dimension $\ge 2$ is called {\em flexible} if its {\em  special
automorphism group} $\SAut(X)$ acts transitively on the smooth locus $X_{\rm
reg}$ \cite{AKZ}. Recall that
$\SAut(X)$ is the subgroup of the
automorphism group $\Aut(X)$ generated by all one-parameter
unipotent subgroups \cite{AKZ}. Given a normal, flexible, affine variety $X$ and
a closed subvariety $Y$ in $X$ of codimension at least $2$,
we show that the pointwise stabilizer subgroup of $Y$ in
the group $\SAut(X)$ acts infinitely transitively on the
complement $X\backslash Y$, that is, $m$-transitively for any $m\ge
1$. More generally we show such a result for any quasi-affine variety $X$ and codimension $\ge 2$ subset $Y$ of $X$.

In the particular case of $X=\AA^n$, $n\ge 2$, this yields a Gromov--Winkelmann Theorem
 \cite{Gr1}, \cite{Wi}.
\end{abstract}
\date{\today}
\maketitle

\thanks{
{\renewcommand{\thefootnote}{} \footnotetext{ 2010
\textit{Mathematics Subject Classification:}
14R20,\,32M17.\mbox{\hspace{11pt}}}}

\vfuzz=2pt
\thanks{}

{\footnotesize \tableofcontents}

\section*{Introduction}

Throughout the paper $X$ will be an algebraic
variety of dimension $\ge 2$ over an algebraically closed field $\kk$ of characteristic 0.
The {\em special automorphism group}
$\SAut(X)$ of such a variety $X$ is the subgroup of the full automorphism group
$\Aut(X)$ generated by all one-parameter unipotent subgroups of
$\Aut(X)$.\footnote{I.e. by subgroups isomorphic to $\G_a$. By abuse of language we do
not distinguish between one-parameter
unipotent subgroups of the group $\Aut(X)$ and effective
$\G_a$-actions on $X$. } Let
$\cU(X)$ denote the set of all these subgroups. A quasi-affine variety $X$ is
called {\em flexible}, if the tangent space $T_xX$ in
any smooth point $x\in X_{\rm reg}$ is spanned by the tangent
vectors at $x$ to the orbits $U.x$, where $U$ runs over $\cU(X)$.

If $X$ is affine then this amounts to the notion of flexibility as introduced in \cite{AKZ, AFKKZ}.
For such varieties the flexibility is equivalent to the transitivity,
and even to infinite transitivity of the group
$\SAut(X)$ acting on the smooth locus
$X_\reg$ of $X$ (see \cite[Theorem 0.1]{AFKKZ}).
(We say that a group action is {\em infinitely transitive} if it is $m$-transitive for any $m\ge 1$.)
These characterizations of flexibility can be extended to any quasi-affine variety
(see  Remarks \ref{1.7} and Theorem \ref{1.10} in Sect.\ 1).

It is worthwhile mentioning that the class of flexible varieties
is rather wide. It includes in particular
\bnum[-]
\item homogeneous spaces of semi-simple groups (and even homogeneous
spaces of extensions of semi-simple groups by unipotent radicals);

\item non-degenerate toric varieties (i.e.\ toric varieties without nonconstant
invertible regular functions);

\item cones over flag varieties and anti-canonical cones over Del Pezzo surfaces of degree at least 4;

\item normal hypersurfaces of the form $uv=p (\bar x )$ in $\C^{n+2}_{u,v, {\bar x}}$;

\item homogeneous Gizatullin  surfaces;
\enum
see \cite{AKZ}, \cite{AFKKZ}, \cite{Kov}.
If on a quasi-affine variety $X$ the group $\SAut(X)$ has an open
orbit, then this open orbit is a flexible quasi-affine variety. A
normal quasi-affine variety $X$ is flexible if and only if so is
$X_\reg$. In its simplest form the main result of this paper is
the following theorem; see Sect.\ 1 for  generalizations and
refinements.

\bthm\label{mthm}
Let  $X$ be a smooth  quasi-affine
variety  of dimension $\ge 2$ and $Y\subseteq X$ a closed subscheme of codimension $\ge 2$. If $X$
is flexible then so is  $X \backslash Y$.
\ethm


That is, if $\SAut(X)$ acts transitively on $X$ then $\SAut(X\backslash Y)$
acts transitively on $X\backslash Y$.
We note that in the setup of the Theorem  any action of a unipotent group on $X\backslash Y$
extends to an action on $X$ preserving $Y$; see Proposition \ref{1.7a} for a more general statement.
Moreover, our main result (see Theorem \ref{mthm1}) yields that the pointwise stabilizer $\SAut_Y(X)$
acts transitively on $X\backslash Y$.
This answers in affirmative a question posed in
\cite[4.22(2)]{AFKKZ}. Partial results in this direction were
obtained in Theorem 2.5 and Proposition 4.19 in \cite{AFKKZ}, see
also Proposition \ref{1.10} below.
Let us note that Theorem \ref{mthm} does not hold for subsets $Y$ of $X$ of codimension 1,
in general; see \cite[Proposition 4.13]{AFKKZ}. In this sense the result above is optimal.

For an affine space $X=\AA^n$, $n\ge 2$,
the flexibility of $X\backslash Y$
was first observed by M.\ Gromov in \cite[\S 2.1.5, p.\ 72,
Exercise (b$'$)]{Gr1}, cf.\ also 4.6(b) and 5.3(c) in \cite{Gr2}.
The transitivity of $\SAut_Y(X)$ in $X\backslash Y$ was proven in
this particular case by J.\ Winkelmann \cite[\S 2, Proposition
1]{Wi}.

The paper is organized as follows. In Section \ref{sec1} we recall
some useful facts from \cite{AFKKZ} and formulate, after
introducing necessary definitions, a stronger version of Theorem
\ref{mthm}, see Theorem \ref{mthm1}.  As an important ingredient of the proof we show that for
any flexible variety $X$ one can find a subgroup 
of $\SAut(X)$ acting with an open orbit on $X$, which is
generated by two locally nilpotent derivations $\delta_0, \delta_1$ along with their replicas $f_0\delta_0$, $f_1\delta_1$, where $f_0\in \ker \delta_0$ and $f_1\in\ker\delta_1$; see Proposition \ref{1.14}. In
Sections \ref{sec2} and \ref{sec3} we prepare the setup for the
proof of Theorem \ref{mthm1}. 

The proof is then contained in Section \ref{sec4}. 
It should be possible, after reading Section
\ref{sec1}, to go directly to  Section \ref{sec4}
addressing results in Sections \ref{sec2} and \ref{sec3} when necessary.

 Let us sketch the scheme of the proof of Theorem \ref{mthm}.
By a result in \cite{AFKKZ} the pointwise stabilizer $\SAut_Y(X)$ of $Y$ in
$\SAut(X)$ has an open orbit, say, $O$ in $X$. We consider a
completion $\bar X$ of $X$ compatible with partial quotients by
the two $\G_a$-subgroups $U_0=\exp(\kk\delta_0)$ and 
$U_1=\exp(\kk\delta_1)$, where $\delta_0$ and $\delta_1$ are as in Proposition
\ref{1.14}.  These quotients define on $\bar X$ two
$\PP^1$-fibrations $\bar\rho_0,\bar\rho_1$ with privileged
sections $D_0,D_1$, which lie on the boundary of $X$ in $\bar X$. Acting with
a suitable replica of $U_0$ one can move the part of the boundary $\partial Y\cap D_1$ to a fixed proper subset of $D_1$, and symmetrically for $U_1$ and $\partial
Y\cap D_0$, see Proposition \ref{prop-one}. Up to a controllable (and so negligible) proper subset
of $D_0\cup D_1$, this property is preserved when we iterate
subsequently actions by  suitable replicas of $U_0$ and $U_1$,  see Proposition \ref{prop-two}.
Using the transitivity property of the subgroup
$H\subseteq\SAut(X)$ generated by $U_0,U_1$ and their replicas, we 
can move a given codimension $\ge 2$ subset $Y$ as in Theorem \ref{mthm}  and,  simultaneously, a given point $x\in X\backslash Y$ to a
generic fiber, say, $F$ of the $\PP^1$-fibration $\bar\rho_0$ so
that $F$ does not meet $\p Y\cap D_0$. Using the Transversality
Theorem from \cite{AFKKZ} we can achieve that $F$ does not meet $Y$ hence in total $F$ and 
$\bar Y$ are disjoint. This enables us to find a $U_0$-invariant function
$f\in\cO_X(X)$, which vanishes on $Y$ and not in $x$. The corresponding replica $U_0'$ of $U_0$ fixes $Y$ and moves $x$ along $F$. Since the
fiber $F$ is generic it meets the open orbit $O$ of $\SAut_Y(X)$, hence so
does  $U_0'.x$. Thus $x$ belongs to $O$, and so $O=X\backslash Y$, as stated.  

In order to prove Propositions \ref{prop-one} and \ref{prop-two} we develop in  Sections \ref{sec2} and \ref{sec3} a machinery, which allows to reduce the proof to the model case of a standard birational transformation of a ruled surface induced by a $\G_a$-action. This reduction is the most lengthly part of the proof.

We thank M.\ Gizatullin for his interest in our work and in particular for his suggestion to treat in Theorem \ref{mthm1} also non-reduced subschemes $Y$ of $X$.

\section{Main theorem}\label{sec1}

\subsection{Basic notions and the main result} We let $\AA^n=\AA^n_\kk$
and $\G_a=\G_a(\kk)$.
In the sequel $X$ denotes a quasi-affine variety over $\kk$.
Thus $X$ can be embedded into an affine variety $X'=\Spec B$ as an open subset.
We let $A=\cO_X(X)$ so that $B$ is a finitely generated $\kk$-subalgebra of $A$. The embedding
$X\emb X'$ factors as $X\to \Spec A\to \Spec B$. Furthermore $X\hto \Spec A$ is an open embedding.
We note that $A$ is in general not a finitely generated algebra over $\kk$.

\blem\la{1.1}
With the notation
as above the following hold.
\bnum[(a)]
\item Every action of an algebraic group on $X$ extends in a canonical way to $\Spec A$.

\item
Every subgroup $U\in \cU(X)$ with infinitesimal generator $\de$ yields a locally nilpotent $\kk$-derivation on $A$.
\enum
\elem

\bproof
(a) is standard, and (b) is a consequence of (a).
\eproof

Let us recall some notions and useful facts from \cite{AFKKZ}.
Given a subgroup $U\in\cU(X)$ we let $\de$
denote an infinitesimal generator of $U$; the latter is uniquely determined up to a nonzero constant factor.
Thus $\de$ is a
locally nilpotent derivation  of the algebra $A=\cO_X(X)$ such that
$U=\exp(\kk \de)$. Geometrically $\de$ can be viewed as a complete
vector field on $X$
with phase flow $u_t=\exp(t\de)$, $t\in \kk$. The tangent vector at the
point $x\in X$ given by this vector field
is denoted $\de_x$.

\blem\la{1.2}
Let $Y$ be a closed (not necessarily reduced)
subscheme of the quasi-affine variety $X$ with ideal sheaf
$\cI\subseteq \cO_X$, and consider the ideal  of global sections $I=\cI(X)\subseteq A=\cO_X(X)$. Given
$U\in \cU(X)$ with an infinitesimal generator $\de$ the following hold.
\bnum[(a)]
\item $\de(A)\subseteq I$  if and only if $u|Y=\id_Y$ for any $u\in U$.

\item $\de(I)\subseteq I$  if and only if $u.Y\subseteq Y$ for any $u\in U$.\footnote{In the terminology
of \cite[p.\ 10]{Fr} this means that $I$ is an integral ideal.}
\enum
\elem

Let us fix the following notation.

\bnota\la{1.3}\label{nota1.3}
(a) Let as before $X$ be a quasi-affine variety and $A=\cO_X(X)$ be its ring of regular functions.
If $\fa\subseteq A$ is the ideal of the complement $\Spec (A)\backslash X$, then the set of nonzero locally
nilpotent derivations $\de$ of $A$ with $\de(\fa)\subseteq \fa$ is denoted by
$$
\LND(X).
$$
In view of Lemmas \ref{1.1} and \ref{1.2}(b) any element
$\de\in \LND(X)$ gives rise to a one-parameter subgroup
$U=\exp(\kk \de)$ in $\cU(X)$ and vice versa.

(b)
In order to deal with quasi-affine varieties we choose a 
 $\kk$-subalgebra $\Lambda$ of $A$ such that
the induced map $X\to \Spec \Lambda$ is an open embedding.
Letting $\fb$ be the ideal of the complement $\Spec(\Lambda)\backslash X$ we let $\LND_\Lambda(X)$ denote the set
of all locally nilpotent derivations $\de$ on
$\Lambda$ with $\de(\fb)\subseteq \fb$. Every such derivation induces as before a one-parameter subgroup
$U\in \cU(X)$ and consequently extends to an element in $\LND(X)$.
Thus $\LND_\Lambda(X)$ can be considered as a subset of $\LND(X)$.

(c) Given a collection $\cN\subseteq \LND_\Lambda(X)$ of nonzero
locally nilpotent derivations we let  $G=G_\cN=\langle\cN\rangle$
be the subgroup of the group $\SAut(X)$
generated by the corresponding one-parameter unipotent subgroups $U=\exp(\kk\de)$, $\de\in \cN$.
\enota

\brems\label{emph} 1.
We emphasize that the subring $\Lambda$ of $A$ is not supposed to be finitely generated over $\kk$
so that the choice $\Lambda=A$ is also possible.
In other words, we consider $X$ as
an open subset of an affine $\kk$-scheme $\Spec \Lambda$,
which is not necessarily an algebraic variety, in contrast with \cite{AFKKZ};
see also
Remark \ref{1.7} below.

2. We observe as well that the $G$-action on $X$ as in \ref{nota1.3}(c)
extends to a $G$-action on the affine scheme $\Spec\Lambda$.\erems

Let us recall some notation and standard facts.

\bsit \la{1.4}
(1) Given a group $G=G_\cN$ as before, the set of all one-parameter
unipotent subgroups of $G$ will be denoted by $\cU(G)$, and the
set of all nonzero locally nilpotent derivations on $\Lambda$ generating
one-parameter subgroups of $G$ by $\LND_\Lambda(G)$ or simply $\LND(G)$.

(2) A {\em $\Lambda$-replica} of a
subgroup $U=\exp(\kk\de)\in\cU(G)$ is a subgroup $U_f=\exp(\kk
f\de)\in\LND_\Lambda(G)$, where $f\in \Lambda$ is in the kernel of $\de$ (\cite{AFKKZ}).

(3) We say that $\cN$
is {\em $\Lambda$-saturated} if $\cN$ is closed under conjugation by
elements in $G$ and taking $\Lambda$-replicas i.e.,
$$
f\de\in\cN\qquad\forall\de\in \cN
\quad\text{and}\quad\forall f\in\ker_\Lambda\de\,.
$$ Hereafter $\Lambda$ will be  fixed, hence in
most cases we omit the symbol $\Lambda$ and say simply `replica' or `saturated'.

(4) A point $x\in X$ is
called {\em $G$-flexible} if $T_xX=\Span(\cN(x))$, where $\cN(x)$ denotes the set of tangent vectors $\de_x$
with $\de\in \cN$. We say that $X$
is {\em $G$-flexible} if $X_{\rm reg}$ consists of $G$-flexible
points.

(5) Given a (not necessarily reduced)  closed subscheme $Y$ in $X$ we let $G_{\cN,Y}$ denote
the subgroup of $G$ generated by all replicas $f\de$ in $\cN$ vanishing on $Y$ in the ideal theoretic sense,
see Lemma \ref{1.2}(a). Therefore
$G_{\cN,Y}\subseteq G_Y$, where $G_Y=\{g\in G:g|Y=\id_Y\}$ stands for the
`pointwise' stabilizer of $Y$ in $G$ in the scheme theoretic sense.
\esit

The following result is our main theorem.

\bthm\label{mthm1} Let $X$ be a  quasi-affine
variety of dimension $\ge 2$ and  $X\hookrightarrow\Spec\Lambda$ be an open embedding into an affine $\kk$-scheme, see \ref{nota1.3}(b).
Let $G=\langle\cN\rangle$
be a subgroup of the group $\SAut(X)$ generated by
a $\Lambda$-saturated set $\cN$ of locally nilpotent
derivations as in \ref{1.4}. Suppose that $X$ is $G$-flexible.
If $Y$ is a closed (possibly non-reduced\footnote{The authors are grateful to M.~Gizatullin for the suggestion
to take also into account non-reduced subschemas $Y$ of $X$.}) subscheme
of $X$ of codimension $\ge 2$, then  the
complement $X \backslash Y$ is $G_{\cN,Y}$-flexible.
\ethm

In the case of a smooth variety $X$ applying Theorem \ref{mthm1} to the group $G=\SAut(X)$ we get
Theorem \ref{mthm}  from  the Introduction.

\brems\la{1.7}
1.  Since $G\subseteq\Aut (\Spec\Lambda)$ the variety $X$ satisfies the requirements of Theorem \ref{mthm1} whenever so does its ($G$-stable) regular locus $X_{\rm reg}$. Therefore
it suffices to prove Theorem \ref{mthm1} under the assumption that $X$ is smooth.
This explains the necessity to fix a subring $\Lambda \subseteq A$ as in \ref{nota1.3}(b). Indeed, $A$ can be properly contained in $A'=\cO_{X_{\rm reg}}(X_{\rm reg})$. If instead of fixing $\Lambda$ we consider always LND's and their replicas with respect to the ring $A=\cO_X(X)$, then an $A'$-replica
is possibly not an $A$-replica and so the notion of saturated set
of derivations could change when passing from $X$
to $X_{\rm reg}$.

2. The viewpoint of the paper \cite{AFKKZ} is slightly different as it deals
with open subsets $X$ of affine algebraic varieties $Z=\Spec B$,
and with subgroups $G$ of $\SAut (Z)$ stabilizing  $X$.
It might happen in principle that although $\Aut(X)$
acts transitively on $X_{\rm reg}$ there is no subgroup $G$ of $\Aut(Z)$
acting transitively on $X_{\rm reg}$, whatever is the choice of an embedding of $X$ into an affine variety $Z$;
cf.\ Question \ref{1.9}
below. Thus {\em a priori}  our viewpoint here is more general. 

3. Working with quasi-affine varieties has yet another advantage: given a subgroup $G\subseteq \SAut(X)$,
in the subsequent proofs we may at any step replace $X$ by an open orbit of $G$. This considerably
simplifies our notation.
\erems

 It is worthwhile to note that if $X$  as in Theorem \ref{mthm1} is normal then the group
$\SAut(X\backslash Y)$ is in a natural way a subgroup of $\SAut(X)$.
This is a consequence of the following proposition.

\bprop\la{1.7a}
Let $X$ be a normal quasi-affine variety and $Y\subseteq X$ a subset of codimension $\ge 2$.
Then every $\G_a$-action on $X\backslash Y$
extends to a $\G_a$-action on $X$ that stabilizes $Y$.
\eprop

\bproof
A $\G_a$-action on $X\backslash Y$ corresponds to a locally nilpotent derivation on $A=\cO_X(X\backslash Y)$
 such that the ideal, say, $\fc$ of the complement $Z\backslash (X\backslash Y)$ is stabilized by $\delta$,
where $Z=\Spec A$.  Because of $\codim_XY\ge 2$
the $\kk$-algebras $A$ and $\cO_X(X)$ coincide.
Consider the ideal $\fa\subseteq A$ of the complement $X^c=Z\backslash X$  and the ideal
$\fb\subseteq A$ of the closure $\bar Y$ so that $\fc=\fa\cap\fb$ is the ideal of the complement of
$X\backslash Y$ in $Z$. We have to show that $\fa$ is stabilized by $\delta$.

This is easy in the case that $A$ is finitely generated, thus $Z$ is an affine algebraic variety.
Indeed, if $U$ stabilizes $X^c\cup \bY$
then it stabilizes all irreducible components of that set 
(see e.g.\ \cite[Proposition 1.14(b)]{Fr}),  thus also $X^c$ and $\bY$ and
consequently their respective ideals.

In the general case, by Lemma \ref{1.8} below $A$ is a direct limit of its $\p$-stable
finitely generated subalgebras $A_i$ such that $X$ embeds as an open subset into $\Spec A_i$.
Applying the first case to every $A_i$ the result follows easily.
\eproof

The following fact is an easy consequence of the Lemma of Cartier \cite[Chapt.\ I, \S 1]{Mu}.

\blem\la{1.8}
Given $\de\in \LND_\Lambda(X)$ and  a finite dimensional $\kk$-subspace $E\subseteq \Lambda$ there is
a finitely generated $\de$-stable $\kk$-subalgebra $\Lambda'\subseteq \Lambda$ containing $E$ such that $X$
embeds as on open subset of the affine variety $\Spec \Lambda'$.
\elem


Since $X$ is quasi-affine there is a finitely generated subalgebra
$C$ of $B$ such that $X$ embeds as an open subset in $\Spec C$. We may suppose that $E$ contains a finite set of generators of $C$.
Since $\p$ is locally nilpotent, the set $E'=\bigcup_{i\ge 1}\p^i(E)$ is finite. Since it  is also $\p$-stable,
it generates a subalgebra $\Lambda'$ of $C$ with  the desired properties.

We do not know whether this result remains true for any finite collection of locally nilpotent derivations.
More precisely:
\smallskip

\bques\la{1.9}
Suppose that $\cN\subseteq \LND(X)$ is a finite subset. Does there exist a finitely generated
$\cN$-stable $\kk$-subalgebra $\Lambda'$ of $A=\cO_X(X)$ such that $X$ embeds into $\Spec \Lambda'$
as an open subset?
\eques

%
%
%
%
%
%
%

\subsection{Transitivity versus flexibility on quasi-affine varieties}

Let $X=\Spec A$ be an affine variety. By the main result in \cite{AFKKZ} the flexibility of $X$
is equivalent to the transitivity of $\SAut(X)$ on $X_\reg$, which in turn is equivalent to infinite
transitivity. In the sequel we need this and related facts in the more general setting
of quasi-affine varieties.

We will state the necessary results in the  generality that we need below. The proofs
in \cite{AFKKZ} can be carried over to our more general quasi-affine setup without any difficulty.
Let us start with the main result of \cite{AFKKZ}, see 1.11 and 2.2 in {\em loc.cit.}

\bthm\label{1.10}
Let $X$ be a smooth, quasi-affine  variety of dimension $\ge 2$,
and let $G=\lan\cN\ran$ be a subgroup of $\SAut(X)$ generated by
a $\Lambda$-saturated set $\cN\subseteq \LND_\Lambda(X)$ as in Notation \ref{1.3} and \ref{1.4}.
Then the following are equivalent.

(i) $X$ is $G_\cN$-flexible.

(ii) $G_{\cN}$ acts transitively on $X$.

(iii) $G_\cN$ acts infinitely transitively on $X$.
\ethm

In the proof of Theorem \ref{mthm1}
we use the following auxiliary results. They are established
in  2.5, 4.19, and  4.2  in \cite{AFKKZ} in the case of  affine schemes $X$ and reduced subvarieties $Y$ of $X$.
The proofs given there carry immediately over to our more general situation.

\bprop\label{1.11}
Let $X$ and $G_\cN$ be as in Theorem \ref{1.10}, and let $Y$ be a closed subscheme of $X$.
If $X$ is $G_\cN$-flexible\footnote{Equivalently, if $G_\cN$ acts transitively on $X$.} then the following hold.

(1) The group  $G_{\cN,Y}$
acts on $X\backslash Y$ with a dense open orbit, say, $O_{Y}$, which
consists of all  $G_{\cN,Y}$-flexible points of $X\backslash Y$.
Consequently, the   $G_{\cN,Y}$-action on $O_{Y}$ is
infinitely transitive.

(2) If $Y$ is finite then $O_{Y}=X\backslash Y$.

(3) If $x\in X$ then
the image of the tangent representation $G_{\cN,x}\to\GL(T_xX)$
given by the differential coincides with the special linear group $\SL(T_xX)$.
\eprop

Finally we need the following interpolation result, see \cite[Theorem 4.14 and Remark 4.16]{AFKKZ}.

\bprop\label{1.12}
Let $X$ and $G_\cN$ be as in Theorem \ref{1.10}. If $G$ acts transitively on $X$
then for any
finite subset $Z\subseteq X$
there exists an automorphism $g \in G$ with $g(x)=x$ for $x\in Z$  and prescribed tangent map $d_xg\in \SL(T_xX)$ at the points $x\in Z$.\footnote{In fact this proposition holds more generally for any finite collection of
$m$-jets provided these jets fix the corresponding points and
preserve local volume forms on $X$ at these points;  see \cite[Remark 4.16]{AFKKZ}.}
\eprop

\subsection{ Generation of subgroups by LND's. }
Let as before $X$ be a quasi-affine algebraic variety of dimension $n\ge 2$
equipped with an open
embedding into an affine $\kk$-scheme $\Spec\Lambda$, where $\Lambda\subseteq\cO_X(X)$.
Given a set of locally nilpotent derivations $\cN\subseteq
\LND_\Lambda(X)$ we enrich it by adding all the $\Lambda$-replicas
of derivations in $\cN$. Letting $\tilde \cN$ be this enlarged set we consider
the subgroup $\lla\cN\rra:=\langle\tilde\cN\rangle$ of the group
$\Aut(X)$ generated by $\tilde\cN$.

In this section we prove the following result.

\bprop\label{1.14} Let $G=\langle\cN\rangle\subseteq\SAut(X)$
be a subgroup generated by a
$\Lambda$-saturated set $\cN$ of locally nilpotent
derivations. Suppose that $G$ acts transitively on $X$.
Then for any locally nilpotent derivation
$\de_0\in \cN$ one
can find another one
 $\de_1\in \cN$ such that the subgroup
 \be\label{3.2a} H=\langle\langle \de_0,\de_1\rangle\rangle
\ee
generated by $\delta_0$, $\delta_1$ and all their replicas acts with an open orbit on $X$.
\eprop

To deduce this result let us recall a few facts.
Let $U$ be a one-parameter unipotent subgroup
with an  infinitesimal generator $\de\in \LND_\Lambda(X)$ (see Notation \ref{1.3}).
By assumption $X$ is contained as an open subset in $\Spec \Lambda$ and by Lemma \ref{1.8} even in
$\Spec \Lambda'$ for some $\de$-stable finitely generated subalgebra $\Lambda'$ of $\Lambda$.
By the Rosenlicht Theorem (see \cite[Theorem
2.3]{PV})
one can find a finite set of $U$-invariant functions $f_1,\ldots,f_m\in \Lambda'^U$,
which separate general $U$-orbits. Let $B$ be the integral closure of the finitely generated $\kk$-algebra
$\kk[f_1,\ldots, f_m]$. It is a standard result that $B$ is again finitely generated,
see e.g.\ \cite[Theorem 4.14]{Ei}.

\bdefi\label{partial quotient}
The  normal affine variety $Q_U=\Spec B$ will be called a
{\em partial quotient} of $X$ by  $U$. In general it depends on the choice of the functions
$f_1,\ldots, f_m$.\footnote{Alternatively, one could  use the {\em Winkelmann quotient} \cite{Wi2}.
This quasi-affine quotient is canonically defined, but has the disadvantage to be non-affine, in general.}
The  inclusion $B\hookrightarrow \cO_X(X)$
defines a dominant morphism $\rho_U:X\to Q_U$ such that the general fibers of $\rho_U$
are general orbits of $U$.
\edefi

\bproof[Proof of Proposition \ref{1.14}]
Let as before $\rho_0:X\to Q_0$ be a partial quotient of $X$ by
$U^0$, where $\dim Q_0=n-1$. Since $n\ge 2$ there exists
$\sigma\in\cN$ such that $\ker\sigma\neq\ker\de_0$ and so $U^0$ and $U=\exp(\C\sigma)$
 have different general orbits.
We can choose $x\in X$ such that the tangent vector  $\de_{0,x}$ of $\de_0$ at $x$ is
nonzero, hence $\dim
U^0.x=1$. Choosing $x$ in an appropriate way there are points $x_1,\ldots,x_{n-1}$ on
the orbit $U^0.x$ such that the vectors
$v_i'=\sigma_{x_i}\in T_{x_i}X$ are all nonzero. Letting $q=\rho_0(x)\in Q_0$ we fix for
each $i=1\ldots, n-1$ a tangent vector $v_i\in T_{x_i}X$ in such a way
that the vectors $d\rho_0(v_i)\in T_qQ_0$, $i=1,\ldots,n-1$, generate the tangent space
$T_qQ_0$ to $Q_0$ at $q$.
For every $i=1,\ldots,n-1$ we can choose a $1$-jet
of a local automorphism at the point $x_i$ that fixes $x_i$ and
sends $v_i'$ to $v_i$. This amounts to choosing
$\alpha_i\in\SL(T_{x_i}X)$ such that $\alpha_i(v_i')=v_i$.
According to Proposition \ref{1.12} one can
interpolate these jets by an automorphism, say, $\alpha\in G$ such
that $\alpha(x_i)=x_i$ and $d\alpha(v_i')=v_i$ for $
i=1,\ldots,n-1$. Replacing $U$ by
$$
U^1=\alpha\circ U\circ \alpha^{-1}=\exp(\C\de_1)\in\cU(G)\,.
$$
we obtain a one-parameter unipotent subgroup with tangent vector $v_i$ at $x_i$, $i=1,\ldots, n-1$.
We claim that the locally nilpotent derivation
$\de_1$ satisfies our requirement. Indeed, $\de_1\in\cN$ since
$\cN$ is saturated  and so, in particular,  is closed under conjugation in $G$.
Consider the conjugated  one-parameter subgroups
$$
U^1_i=\alpha_i^{-1}\circ U^1\circ \alpha_i=\exp(\C\sigma_i)\in\cU(H),\quad
i=1,\ldots,n-1\,,
$$
where $\alpha_i\in U^0$ is an element which maps $x$ to $x_i$. Here $H$ is as in (\ref{3.2a}) and
$\sigma_i$ is a conjugate of $\de_1$ under the action of $H$ for $ i=1,\ldots,n-1$. For any
$i$ in this range the vector $u_i=d\alpha_i(v_i)$ is tangent to the
orbit $U^1_i.x$ at the point $x\in X$. Furthermore, the
vectors $d\rho_0(u_i)=d\rho_0(v_i)\in T_q  Q_0$, $i=1,\ldots,n-1$
still generate $T_qQ_0$. Hence the vectors
$$
u_0=\de_{0,x},\,\,u_1=\sigma_{1,x},\,\,\ldots,\,\,u_{n-1}
=\sigma_{n-1,x}\in T_xX
$$
span $T_xX$ as well.
Consequently, $x$ is an $H$-flexible point and so the $H$-orbit
$H.x$ is open and dense in $X$ (see \cite[Corollary 1.11(a)]{AFKKZ}).
\eproof

\section{$m$-blowups, tangency, and $m$-contractions}\la{sec2}
This section is technical; we use its results and notions (see especially Definitions \ref{x.5} and \ref{x.8} and Proposition \ref{x.13})
in the proof of Proposition \ref{prop-one} in the next section. 

\bsit\la{x.1}
In the  sequel we  deal with rational maps $\bdi g:X&\rDotsto& Y\edi$ which
fit into a diagram
\bdi
&& \hat X\\
&\ldTo<{h}&&\rdTo>{g'}\\
X&&\rDotsto^g && Y
\edi
where $h$ is a sequence of blowups and $g'$ is a proper morphism. This somewhat restricted class of
rational maps is suitable for our purposes. Given subsets $A\subseteq X$ and $B\subseteq Y$ we let
$$
g(A)=g'(h^{-1}(A))\and g^{-1}(B)=h(g'^{-1}(B))
$$
denote the total image and preimage, respectively.\footnote{These notions should be treated with
caution, because they are not compatible with composition of rational maps.} Since any two resolutions
of the indeterminacy set are dominated by a third one, the total image
and the total preimage are well defined.
\esit

\subsection{$m$-blowups and tangency}
In the next Definition we introduce a setup which is used repeatedly in this and the next section.

\bdefi\label{x.2}
Let $X$ be an  algebraic variety and $C$, $D$ be divisors in $X$, which are Cartier near $C\cap D$.
The {\em $m$-blowup $\sigma_m:X_m\to X$ of $D$ along $C$} is defined recursively as follows.
With $X_0=X$ we let $X_1$ be the blowup of $X$ along the subscheme $C\cap D$.
 If $X_{m-1}$ is already defined for some $m\ge 2$,
 then we let $X_{m}\to X_{m-1}$ be the blowup along $D^{(m-1)}\cap
E_{m-1}$, where $D^{(m-1)}$ is the proper transform of $D$ in $X_{m-1}$ and $E_{m-1}$ the exceptional
set  of the previous blowup $X_{m-1}\to X_{m-2}$.

In the following we call the proper transforms
$$
E'_1,\ldots, E_m'\subseteq X'= X_m
$$
of the exceptional sets $E_i$ of $X_i\to X_{i-1}$ the {\em exceptional sets of the $m$-blowup} of $D$
along $C$. The proper transforms of $C$ and $D$ will always be denoted $C',D'$, respectively.
\edefi

\bexa\la{x.3}
Suppose that $S$ is a complete smooth surface and $C\cap D=\{p\}$, where the intersection
is transversal. Then the dual graph of
$C'\cup E_1'\cup\ldots\cup E_m'\cup D'$ is a linear chain:
\be\label{eqegraph} \cou{C^2-1}{C'}\llin
\cou{-2}{E'_1}\lin\ldots\lin
\cou{-1}{E'_m}\vlin{14}\cou{D^2-m}{D'}
\qquad.\ee
\eexa

Let us consider next the effect of an $m$-blowup as in Definition \ref{x.2} on the boundary of
a closed subset of $X$.

\bprop\label{x.4}
We keep the notation and assumptions as in Definition \ref{x.2}.
Given a closed subset $Y\subseteq X$ we let $Y'$ denote its proper transform in $X'$ and $\p Y'$
its boundary $\p Y'= Y'\cap \sigma_m^{-1}(C\cup D)$.
Then with $P= \ol{Y\cap D\backslash C}$, for $m\gg 0$
$$
\p Y'\subseteq E_1'\cup\ldots\cup E_{m-1}' \cup \sigma_m^{-1}(P)\,.
$$
\eprop

\bproof
The assertion is local around points in $C\cap D\backslash P$. Thus we may assume
that $P=\emptyset$, $X=\Spec A$ is affine, and that $D=V(x)$, $C=V(y)$
with functions $x,y\in A$. The subset
$$
U'=X'\backslash \bigcup_{i=0}^{m-1}  E'_i
$$
of $X'$ is affine with coordinate ring
$$
A'= A[u]\,,\quad\mbox{where }u=x/y^m,
$$
cf.\ Lemma \ref{x.10} below for the special case of surfaces. Furthermore
\be\la{aux6a}
U'\cap E'_m=\{y=0\}\and U'\cap D'=\{u=0\}.
\ee
If $I\subseteq A$ is the ideal of $Y$ then $B=A/I$ is the affine coordinate ring of $Y$.
Since $\ol{Y\cap D\backslash C}=\emptyset$ the set $Y\cap D$ is contained in $C\cap D$
and so the localization $(B/xB)_y$ is zero. Hence there exists a  natural number $m$
such that $y^{m-1}\in xB$. In other words, we can find $a\in A$ such that
\be\la{aux6b}
y^{m-1} -a\cdot x \in I \,.
\ee
In the blowup ring $A'$
the ideal $I'$ of $ Y'$ is given by
$$
I'=\{g\in A'\mid \exists k\in \NN: y^kg\in IA'\}.
$$
Since $u=x/y^m$ condition \eqref{aux6b} can be rewritten in the form
$$
y^{m-1}\cdot\left(1- yau\right)\in IA'.
$$
Hence $1- yau\in I'$. This shows that in the affine coordinate ring $B'=A'/I'$ of $U'\cap  Y'$
the residue classes of  $y$ and $u$ are units. In view of \eqref{aux6a} this implies that
$$
U'\cap  Y'\cap  E'_m=\emptyset\and
U'\cap Y'\cap  D'=\emptyset,
$$
which immediately yields the required result.
\eproof

\bdefi\la{x.5}
We say that a closed subset $Y$ of $X$ is {\em at most $m$-tangent} to $D$ along $C$,
if the conclusion of Proposition \ref{x.4} holds with this particular value of $m$.
The subset $N=C\cap \ol{Y\cap D\backslash C}$ of $C\cap D$ will be called the {\em defect set}.
\edefi

We note that if $Y$ is at most $m$-tangent to $C$ along $D$ then it is also at most $m'$-tangent
to $C$ along $D$ for all $m'\ge m$.
The following observation is important.

\blem\la{x.6}
If $\codim_XY\ge 1$ and $Y\backslash D$ is dense in $Y$ then the defect set $N$ is nowhere
dense in $C\cap D$.
\elem

\bproof
If $\codim_XY\ge 1$ then the set $Y\cap D$ has
codimension $\ge 1$ in $D$. Hence its closure
cannot contain any component of $C\cap D$. 
\eproof

\brem\la{x.7}
In the setup of Proposition \ref{x.4} suppose that $(Y_s)_{s\in S}$ is a family of proper
closed subsets of $X$. Then there is a natural $m$  such that $Y_s$ is at most $m$-tangent
to $D$ along $C$ for any $s\in S$.

This follows easily from the fact that the construction of  Proposition \ref{x.4} can be
done al least generically in the given family and that it is then compatible with restriction
to the general fiber. More precisely, one can find an open dense subset $U\subseteq S$
so that all fibers $Y_s$ are at most $m$-tangent to $D$ along $C$ with $m$ independent
of $s\in U$, and with a defect set $N_s=C\cap\ol{Y_s\cap D\backslash C}$.
Restricting the family to $S'=S\backslash U$ and applying induction on $\dim S$,
we may assume that $Y_s$ is at most $m$-tangent to $D$ along $C$ for any $s\in S'$.
Hence the assertion follows.
\erem

\subsection{$m$-contractions}

\bdefi\la{x.8}
Let $C$, $D$ be divisors  on an the algebraic variety $X$,
which are Cartier near $C\cap D$. Consider a  birational map $\bdi g:X&\rDotsto & X\edi$
and a resolution  of the indeterminacy set of $g$ which factors through the $m$-blowup
$\sigma_m:X'=X_m\to X$ of $D$ along $C$, see  Definition \ref{x.2}:
\bdi
&& \hat X\\
&\ldTo<{h_m}&\dTo>h &\rdTo>{g'}\\
X'=X_m& \rTo^{\sigma_m} &X&\rDotsto^g & X
\edi
$g$ is called an {\em $m$-contraction for $C$ along $D$} if the following hold.
\bnum
\item $g$ is biregular in the points of $X\backslash C$;
\item  with $g_m=g\circ\sigma_m$, the total image\footnote{See \ref{x.1}}
$g_m(C'+E_1'+\ldots+ E_{m-1}')$ is  a subset of $D$, where   $E_1', \ldots, E_m'$ are
as in Definition \ref{x.2}.
\enum

\edefi

Clearly, an $m$-contraction for $C$ along $D$ is also an $m'$-contraction for $C$ along
$D$ for any $m'\le m$.
The following  example is important and serves as a model case.

\bnota\la{x.9}
Let $\Gamma=(\Gamma,o)$ be a germ of a smooth affine curve with a uniformizing parameter $u$ such that $u(o)=0$,
and let $d(u)$ denote a  nowhere vanishing function on $\Gamma$.  We consider homogeneous coordinates $(\zeta_1:\zeta_2)$ on $\PP^1$ and an affine coordinate
$v=\zeta_1/\zeta_2$ on $\AA^1=\PP^1\backslash \{(1:0)\}$. The product $S:=\Gamma\times\PP^1$  is a $\PP^1$-fibered surface over
$\Gamma$. Its fiber, say, $C$ over $o\in \Gamma$ and the section
$D=\Gamma\times \{(0:1)\}\subseteq \Gamma\times\AA^1$ can be described in  coordinates by 
$$
C=\{u=0\}\and D=\{v=0\}.
$$
 Let us study the rational map $g_m: S\dashrightarrow  S$, where $m\in\NN$,
given in affine coordinates by
\be\la{hm}
g_m(u,v)=
\left( u,\,\,\frac{u^mv}{d(u) v+u^m} \right)\,.
\ee
Its indeterminacy set consists of the intersection point $C\cap D=\{u=v=0\}$, which will be
denoted by $\bar 0$.
\enota

\blem\la{x.10}
Let
\bdi[small]
& &{S'} & &\\
&\ldTo^{\sigma_m} &  &
\rdTo>{{g'}_m}&  \\
 S &&\rDashto^{g_m} &&  S\\
\edi

\noindent
be the minimal resolution of indeterminacies of $g_m$, where  $\sigma_m$ is a sequence of blowups
and $g'_m$ is a morphism.
Then the total transform of $C+D$ on $S'$ under $\sigma_m$ has weighted dual graph
\medskip

\be\label{aux8a} \cou{-1}{C'}\lin
\cou{-2}{E'_1}\lin\ldots\lin
\cou{\!\!\!\!\!\!\!\!\!\!\!\!-2}{E'_m}\nlin\cshiftup{D'}{-m}
\lin\ldots\lin\cou{-2}{E'_{2m-1}}\lin\cou{-1}{E'_{2m}} \quad,\ee
\noindent
where  $C'$ and $D'$
are the proper transforms of $C$
and $D$, respectively.
The map  $\sigma_m$ contracts the components
$E'_1,\ldots,E'_{2m}$ to the origin
$\bar 0\in S$, while $g'_m$ contracts the curves
$C',E'_1,\ldots,E'_{2m-1}$ to $\bar 0\in S$. Furthermore $g'_m(D')=D$ and $g'_m(E'_{2m})=C$.
\elem

\bproof
Letting $v_0=v$ we define a sequence of coordinates charts $(u,v_i)$
on $S'$, $i=0,\ldots,2m$, so that the
$2m$ blowing-downs over the origin
with exceptional curves $E'_1,\ldots,E'_{2m}$ that constitute the map
$$
\sigma:(u,v_{2m})\mapsto (u,v_{2m-1})\mapsto\ldots
\mapsto(u,v_{1})\mapsto (u,v)\,
$$
can be described by the formulae
\be\la{x.10a}
v_1=v/u,\quad v_2=v_1/u=v/u^2,\quad\ldots,\quad
v_m=v_{m-1}/u=v/u^m\,,
\ee
and
\be\la{x.10b}
v_{m+1}=(1 +d(u) v_m)/u,\quad v_{m+2}=v_{m+1}/u,\quad
\ldots,\quad v_{2m}=v_{2m-1}/u=(1+d(u)v_m)/u^m\,.
\ee
The map $g_m$ can be written in these coordinate charts as
$$
(u,v)\mapsto \left(u,\,\frac{u^mv}{d(u)v+u^m}\right)=
\left(u,\, \frac{u^mv_1}{d(u)v_1+u^{m-1}}\right)=\ldots
$$
$$
\ldots =
\left(u,\, \frac{u^m v_m}{1+d(u) v_m}\right)
=\left(u,\, \frac{d(u) u^mv_{m+1}- u^{m-1}}{v_{m+1}}\right)=\ldots=
\left(u,\,\frac{d(u) u^m v_{2m}-1}{v_{2m}}\right)\,.
$$
Hence the curve $E'_i$ given in the chart $(u,v_i)$ by equation $u=0$ is contracted under
$g'_m$ for every $i=0,\ldots,2m-1$,
while the curve $E'_{2m}$
given by the same equation in the chart $(u,\, v_{2m})$
maps birationally onto the curve $C$ in $S$. Now the assertion
follows.\eproof

An immediate consequence is the following corollary.

\bcor\la{x.11}
The birational map $g_m$ in \eqref{hm} is an $m$-contraction of $C$ along $D$.
\ecor

Let us note that  $g_m$ is not an ($m+1$)-contraction of $C$ along $D$.
This example can be generalized to higher dimensions as follows.

\bnota\la{x.11a}
Instead of a curve $\Gamma$ in \ref{x.9} we consider now a smooth affine algebraic variety $Q$
and a smooth divisor $T\subseteq Q$ given by the equation $\{u=0\}$, where $u\in \cO_Q(Q)$.
The product $X=Q\times \PP^1 $ is $\PP^1$-fibered over $Q$ and contains the divisors
$$
C=T\times \PP^1\and
D=Q\times \{(0:1)\}\subseteq Q\times\AA^1,
$$
where we equip $\PP^1$ with homogeneous coordinates $(\zeta_1:\zeta_2)$. As before
$v=\zeta_1/\zeta_2$ stands for an affine coordinate on $\AA^1=\PP^1\backslash \{(1:0)\}$.
Thus  we have
$$
C=\{u=0\}\and D=\{v=0\}.
$$
\enota

\blem\label{x.11b}
Given a nowhere vanishing function $d(q)$ on $Q$ and $m\in\NN$  the rational map
\be\la{x.11aa}
g_m: X\dashrightarrow  X\,,\quad \mbox{where}\quad
g_m(q,v)=
\left( q,\,\,\frac{u(q)^mv}{d(q) v+u(q)^m} \right)\,,
\ee
is an $m$-contraction of $C$ along $D$.
\elem

\bproof
A resolution
\bdi[small]
& &{X'} & &\\
&\ldTo^{\sigma_m} &  &
\rdTo>{{g'}_m}&  \\
 X &&\rDashto^{g_m} &&  X\\
\edi
of the indeterminacy points of $g_m$ can be obtained (with obvious changes) by the same sequence
of blowups as in the proof of Lemma \ref{x.10}. Letting $v_0=v$ we define a sequence of coordinates
charts $(q,v_i)\in U_i=Q\times \AA^1$
on $X'$, $i=0,\ldots,2m$, so that the
$2m$ blowdowns over $C\cap D$
with exceptional divisors $E'_1,\ldots,E'_{2m}$ that constitute the map
$$
\sigma:(q,v_{2m})\mapsto (q,v_{2m-1})\mapsto\ldots
\mapsto(q,v_{1})\mapsto (q,v)\,
$$
can be described by the formulae in \eqref{x.10a} and \eqref{x.10b}, where $u$ is now the function $u(q)$.
With the same calculation as before the map $g_m$ can be written in these coordinate charts as
$$
(q,v)\mapsto
\left(q,\,\frac{d(q) u(q)^m v_{2m}-1}{v_{2m}}\right)\,.
$$
As in the proof of \ref{x.10} the exceptional set $E'_i$ is given in the chart $U_i$ by the equation $u=0$,
and it is contracted under $g'_m$ to the subset $C\cap D$ for every $i=0,\ldots,2m-1$.
Finally, the exceptional set $E'_{2m}$
given by $\{u=0\}$ in the chart $U_{2m}$
maps under $g_m'$ isomorphically onto the divisor $C$ in $X$. Since
the divisors $C', E_1',\ldots, E_{m-1}'$ in $X'$ are contracted under $g_m'$ to $C\cap D$,  the result follows.
\eproof

Next we show that $m$-contractions are compatible with certain blowups.

\bprop\la{x.12}
Let $X$ be an algebraic variety and $C$, $D$ be connected divisors on $X$, which are
Cartier near $C\cap D$.
Let $\bdi g: X&\rDotsto &X\edi $ be an $m$-contraction of $C$ along $D$ and $p:Z\to X$ be
a modification, which is  an isomorphism over $D\cup (X\backslash C)$. Then  the rational map
$\bdi f: Z&\rDotsto &Z\edi $ induced by $g$ is an $m$-contraction of $C_Z=p^{-1}(C)$ along
$D_Z=p^{-1}(D)\cong D$.
\eprop

\bproof
Let $X_m\to X$ and $Z_m\to Z$ be the $m$-blowups of $X$ and $Z$, respectively. Since $p$ is an
isomorphism in the points near $D$, the exceptional sets $E_1',\ldots, E_m'$
 of $X_m\to X$ can be identified in a natural way with the exceptional sets, say,
$E_{1,Z}',\ldots, E_{m,Z}'$
of $Z_m\to Z$.  Consider the composed rational maps
\bdi
Z'_m &&\rDotsto^{f_m} && Z\and X'_m &&\rDotsto^{g_m} && Z\,.
\edi
and a diagram
\bdi
&& \hat Z\\
&\ldTo<{h_m}&&\rdTo>{f_m'}\\
Z'_m &&\rDotsto^{f_m} && Z\\
\dTo<{p'} &&&&\dTo>p\\
X'_m &&\rDotsto^{g_m} && X
\edi
where $\hat Z$ is a resolution of the indeterminacy locus of $f_m$
and then also of $g_m$. By our assumption the set
$$
(p'\circ h_m)^{-1}(C'\cup E_1'\cup\ldots\cup E'_{m-1})
=h_m^{-1}(C_Z'\cup E_{1,Z}'\cup\ldots\cup E'_{m-1,Z})
$$
 is contracted under $p\circ f_m'$ to a subset of $D$. Since $p$ is an isomorphism near $D$
the latter set is already contracted under $f_m'$ to a subset of $D$. This proves the assertion.
\eproof

Let us now study the effect of an $m$-contraction of $C$ along $D$ on the boundary
of a closed subset $Y$ of $X$.

\bprop\la{x.13}
Let $X$ be an algebraic variety and $C$, $D$ divisors on $X$, which are Cartier
near $C\cap D$.
Assume that $\bdi g:X&\rDotsto& X\edi$ is an $m$-contraction of $C$ along $D$ and
that  $Y\subseteq X$
is a closed subset, which is at most $m$-tangent to $C$ along $D$ with defect set
$N=C\cap \ol{Y\cap D\backslash C}$.
Then the proper image $\hat Y$ of $Y$ under $g$ satisfies
$$
\p\hat Y\subseteq  D\cup g(N)\,,
$$
where $g(N)$ is the total image of $N$ and $\p \hat Y$ denotes the intersection of $\hat Y$
with $D\cup C$.
\eprop

\bproof
Let $\sigma_m:X'=X_m\to X$ be the $m$-blowup of $C$ along $D$ with exceptional sets
$E'_1,\ldots, E'_m$ and consider the composition $\bdi g_m=g\circ\sigma_m:X'&\rDotsto&X\edi$.
We can find a resolution of the indeterminacy locus of $g_m$
\bdi
&& \hat X\\
&\ldTo<{h_m}&&\rdTo>{g_m'}\\
X'&&\rDotsto^{g_m} && X\,.
\edi
Since $Y$ is at most $m$-tangent to $C$ along $D$, the boundary $\p Y'$ of the proper transform $Y'$
of $Y$ in $X'$ satisfies
$$
\p Y'\subseteq E_1'\cup\ldots\cup E_{m-1}'\cup \sigma_m^{-1}(P)\,,
$$
where $P=\ol{Y\cap D\backslash C}$, see Proposition \ref{x.4}.
By condition (2) in Definition \ref{x.8}
$$
h_m^{-1}(C'\cup E_1'\cup\ldots\cup E_{m-1}')
$$
is contracted under $g_m'$ to a subset of $D$. Hence
$$
g_m'(h_m^{-1}(\p Y'))\subseteq D\cup g_m'(h_m^{-1}(\sigma_m^{-1}(P))))
=D\cup g(P).
$$
Since $g_m'$ is proper the set on the right is easily seen to contain $\p \hat Y$, as stated.
\eproof

\section{Replicas as $m$-contractions}\la{sec3}

\bnota\label{5.1}
(a) Let $X$ be a smooth quasi-affine algebraic variety and $G_\cN$ a group of automorphisms on $X$
generated by a set of $\Lambda$-saturated locally nilpotent derivations
$\cN\subseteq \LND_\Lambda(X)$,
see Notation \ref{1.3} and \ref{1.4}. Suppose that $G_\cN$ acts transitively on $X$.

(b) We choose two locally nilpotent derivations $\de$, $\de_0\in \LND_\Lambda(X)$ such that
$$
\ker\de\ne \ker\de_0 .
$$
Let $U$, $U^0$ denote the associated one-parameter subgroups and choose partial quotients
$$
\rho:X\to Q \and \rho_0:X\to Q_0
$$
as introduced  in \ref{partial quotient}.

(c) We can embed $Q$ and $Q_0$ into normal projective  varieties $\bar Q$ and $\bar Q_0$,
respectively. Let $\bX$ be a smooth projective completion of $X$. After blowing up $\bX$ in the boundary $\p X=\bX\backslash X$, if necessary, we may extend $\rho$ and $\rho_0$  to morphisms
\bdi
\bX&\rTo^{\brho_0}& \bQ_0\\
\dTo>{\brho}\\
\bQ
\edi
The general fiber of $\rho$ is an orbit of $U$
isomorphic to $\AA^1$. Clearly
$$
\brho^{-1}(q)\cong \PP^1
$$
for a general point $q\in Q$. Hence there is a unique divisor  $D\subseteq \bar X\backslash X$
which maps  birationally onto $\bar Q$.
Similarly there is a unique divisor $D_0$ in $\bX\backslash X$
mapping birationally
onto $\bQ_0$. Thus both $D$ and $D_0$ are contained in the boundary $\p X=\bX\backslash X$. \enota

The following  observations will be important.

\blem\label{5.2}
\bnum[(1)]
\item Let $\vp\in \ker \de \backslash \ker \de_0$ be a regular function on $X$.
Then $\vp$ is a rational function on $\bX$ with poles at general points of $D_0$.

\item We have
$$
\brho(D_0)\subseteq \bQ\backslash Q
\and
\brho_0(D)\subseteq \bQ_0\backslash Q_0\,.
$$
In particular, $D\ne D_0$.
\enum
\elem

\bproof
(1) Since $D_0\to \bQ_0$ is dominant, an orbit closure $\overline{H_0 .x}$ of a general point $x\in X$
meets $D_0$ at a general point $\bar x\in D_0$. Let us consider
$\vp$ as a rational map
$\bdi \bX&\rDashto & \PP^1\edi$. Since the indeterminacy set of $\vp$ on $\bX$ is of codimension
at least 2, $\vp$ is regular on the orbit closure $\overline{H_0 .x}\cong \PP^1$ for a general $x\in X$.
Since $\vp\not\in \ker \de_0$ this map is not constant on general orbits of $H_0$.
In particular it restricts to a dominant morphism $\varphi:\overline{H_0 .x}\to\PP^1$ such that
$\varphi(\bar x) = \infty$.

(2) It is sufficient to prove the first part. If $\brho(D_0)\cap Q\ne\emptyset$ then a function
$\vp\in \cO(Q) \backslash \ker \de_0$ would be holomorphic in a general point of $D_0$ contradicting  (1).
\eproof

\blem\la{5.3}
After blowing up the boundaries $\p X=\bX\backslash X$ and $\p Q=\bQ\backslash Q$
suitably we can achieve that
\bnum[(a)]
\item $T=\brho(D_0)$ is a divisor in $\bQ$, and
\item $\bX$, $D$ and $D_0$ are smooth.
\enum
\elem

\bproof
(a) By Lemma \ref{5.2}(2) $T$ sits in the boundary of $\bQ$.
According to a theorem of Zariski, see \cite{Za} and \cite[Theorem 1.3]{Kol}, there is a  blowup  $\bar Q'\to\bar Q$ with a center in
$\brho(D_0)$ such that the proper transform of $D_0$ in $\bX_{\bQ'}$ maps onto a divisor in $\bQ'$.
Thus replacing $\bar Q$ by $\bar Q'$ we can achieve that $T$ is a divisor.

Since $X$ is smooth and does not meet $D\cup D_0$, by a suitable blowup of the boundary $\bX\backslash X$
we can achieve that (b) holds.
\eproof

\blem\la{5.3a} \la{}
There is a closed subset $B_0$ of $\bQ$ with $\codim_\bQ B_0\ge 2$   such that the following hold.
\bnum[(a)]
\item $\Sing \bQ\cup \Sing T\subseteq B_0$.
\item $D\to \bQ$ is an isomorphism in the points $D\backslash \brho^{-1}(B_0)$.
\item $\bX\to \bQ$ is flat in the points over $\bQ\backslash B_0$.
\enum
\elem

\bproof 
(a) can be satisfied as $\bQ$ is normal and $T$ is reduced. Since $D\to \bQ$ is a birational map, also (b) can be achieved.

(c)
By the theorem on generic flatness \cite[Theorem 14.4]{Ei}
there is a proper closed subset $E$ in $\bQ$ such that $\brho$
is flat in the points over $\bQ\backslash E$. Applying the theorem on generic flatness again gives that
the restricted map $\brho|E:\brho^{-1}(E)\to E$ is flat over a subset $E\backslash B'$ of $E$, where $B'$
is a nowhere dense closed subset of $E$. Using Corollary 6.9 in \cite{Ei} it follows that $f$ is flat over
the set $\bQ\backslash B''$, where
$$
B''=B'\cup \{s\in E: E \mbox{ is not a Cartier divisor in $\bQ$ at }x \}
$$
Since $\bQ$ is normal this set has codimension $\ge 2$ in $\bQ$.
Adding $B''$ to $B_0$, also (c) is satisfied.
\eproof

The following facts should be well known; in lack of a reference
we provide a brief argument.

\blem\la{5.4}
Let $p:S\to \Gamma$ be a $\PP^1$-fibration of a smooth surface $S$ over a smooth affine curve $\Gamma$
admitting a smooth section $D\subseteq S$ so that $D\cong \Gamma$. Then for any point $t\in\Gamma$
the fiber $F=p^{-1}(t)$ over $t$ is a tree of rational curves. Furthermore the following hold.
\bnum[(a)]
\item  
If $\{x\}=F\cap D$ then $h^0(F,\cO_F(x))=2$ and $H^i(F,\cO_F(x))=0$ for $i\ge 1$.

\item The sheaf $\cO_F(x)$ is generated by its global sections.

\item If $s_0, s_1\in H^0(F,\cO_F(x))$ is a basis, then the map $(s_0:s_1): F\to \PP^1$ is an isomorphism near $x$.
\enum
\elem

\bproof
Blowing down successively $(-1)$-curves in the fibers of $p$ not meeting $D$ we obtain a locally trivial $\PP^1$-bundle
$\cV\to \Gamma$.  The curve $D$ can as well be considered as a section of $\cV\to \Gamma$ and so we have an isomorphism
$\cV\cong\Proj_\Gamma(p_*(\cO_\cV(D)))$. If $S=\cV$ then the assertions (a)-(c) are trivial.
Blowing up subsequently points in the fibers these assertions also follow for $p:S\to \Gamma$. 
\eproof

In what follows we may assume that the conditions (a), (b) in Lemma \ref{5.3} are satisfied.

\blem\la{5.5}
Letting $\bX_q=\brho^{-1}(q)$ and $D_q=D\cap \bX_q$ 
there is a  closed subset $B$ of codimension $\ge 2$ in $\bQ$ such that for $q\in \bQ\backslash B$  the following assertions hold.
\bnum[(a)$_q$]
\item $h^0(\bX_q,\cO_{\bX_q}(D_q))=2$ and $H^i(\bX_q,\cO_{\bX_q}(D_q))=0$ for $i\ge 1$.

\item The sheaf $\cO_{\bX_q}(D_q))$ is generated by its global sections.
\item If $s_0, s_1\in H^0(\bX_q,\cO_{\bX_q}(D_q))$ is a basis, then the map $(s_0:s_1)): \bX_q\to \PP^1$
is an isomorphism near $D_q$.
\item The map $\brho_*(\cO_{\bX}(D))_q\to
H^0(\bX_q,\cO_{\bX_q}(D_q))$ is surjective, and $\brho_*(\cO_{\bX}(D))_q$ is free of rank 2.
\enum
\elem

\bproof
Let $B_0\subseteq \bQ$ be a set as in Lemma \ref{5.3a}.
We choose a proper closed subset $P$ of $\bQ$ such that any fiber over $\bQ\backslash P$ is isomorphic to $\PP^1$.
For any $q\in \bQ\backslash P$ the assertions (a)$_q$-(d)$_q$ follow easily.

Let a curve $\Gamma$ in $\bQ$ be an intersection of $n-1$ general ample divisors in $\bQ$.
Since $\bQ$ is normal and $\codim B_0\ge 2$, $\Gamma$ meets neither $\Sing \bQ$ nor $B$.
By Bertini's theorem both $\Gamma$  and the surface $S=\brho^{-1}(\Gamma)$ are smooth.
The restriction $\brho|S:S\to\PP^1$ is a $\PP^1$-fibration. This $\PP^1$-fibration  admits a section,
namely $D\cap S$. The intersection  $D\cap S$ is smooth in view of Bertini's theorem and Lemma
\ref{5.3}(b). The fiber of $S\to \Gamma$ over $q\in \Gamma\subseteq \bQ$ coincides with $\bX_q$.
By Lemma \ref{5.4} such a fiber $\bX_q$ is a tree of rational curves satisfying (a)$_q$-(c)$_q$. Since
$\Gamma$ meets every component, say, $P_i$ of $P$ of codimension 1 and does not meet $B_0$,
for some $q_i\in P_i\backslash  B_0$ the conditions (a)$_{q_i}$-(c)$_{q_i}$ are satisfied.
By semicontinuity (see\cite[III, 12.8]{Ha}) we obtain the inequalities
$$
h^j(\bX_p,\cO_{\bX_p}(D_p))\le
h^j(\bX_q,\cO_{\bX_q}(D_q))\le h^j(\bX_{q_i},\cO_{\bX_{q_i}}(D_{q_i})), \quad j\ge 0,
$$
where $q\in P_i$ is a point near $q_i$ and $p\in \bQ\backslash P$ is a point near $q$.
Since the outer terms are equal, condition
(a)$_q$ holds for $q$ in some open dense subset $P^o_i$ of $P_i$.

By Grauert's criterion (see \cite[III, 12.9]{Ha}) now also (d)$_q$ is satisfied. Since (b)$_q$ and (c)$_q$
are open conditions on $P_i^o$, which are satisfied for some $q\in P_i^o$, they are satisfied generically
on $P_i$. Now the lemma follows.
\eproof

\bcor \la{5.6}
There is a proper closed subset $B\subseteq \bQ$ containing $\Sing T$ and $\Sing \bQ$ with
$\codim_T(T\cap B)\ge 1$ such that, letting
$$
X^o=\bX\backslash \brho^{-1}(B)\,,
\quad Q^o=\bQ\backslash B\,,\quad  T^o=T\backslash B\and
C=\brho^{-1}(T)\,,
$$
there is a birational morphism
\be\la{5.6a}
\vp: X^o\lto \cX= Q^o\times \PP^1
\ee
compatible with the projection to $Q^o$,
which restricts to a biregular morphism
\be\la{5.6b}
X^o\backslash C\lto \cX\backslash \cC=(Q^o\backslash T)\times \PP^1\,,
\ee
where $\cC=T^o\times \PP^1$. Furthermore $\vp$ is biregular in a neighborhood of $D^o=D\cap X^o$.
\ecor

\bproof
Let $B\subseteq \bQ$ be the subset constructed in Lemma \ref{5.5}. Enlarging it in a suitable
way we may assume that it contains $\Sing T\cup\Sing \bQ$. According to Lemma \ref{5.5}(c)
the sheaf $\cE=\brho_*(\cO_{X^o}(D))$ is locally free of rank 2 on $Q^o$. Thus enlarging $B$
we may suppose that $\brho_*(\cO_{X^o}(D))$ is free. Choose two sections $s_0, s_1$ which
form a basis of this bundle. They provide a morphism
$$
\vp= (\brho, (s_0:s_1)): X^o\to Q^o\times \PP^1.
$$
Restricting to a fiber over $q\in Q^o$, in view of Lemma \ref{5.5}(c)$_q$ this yields an isomorphism
near $D_q$.
Hence $\vp$ is an isomorphism near $D^o$.
Enlarging $B$ further we may also assume that all fibers in $Q^o\backslash T$ are isomorphic to $\PP^1$.
This implies that the restricted morphism \eqref{5.6b} is an isomorphism.
\eproof

\bnota\la{5.7}
Consider the restriction of the locally nilpotent vector field $\de$ to $X^o\cap X$.
The associated action of $U=\exp(\kk\de)$ has no fixed points in this set and  extends to an action on
$X^o\backslash C$, where as before $C=\brho^{-1}(T)$. The fibers of $X^o \backslash C\to Q^o\backslash T$
are preserved under $U$.

Under the isomorphism $X^o\backslash C \simeq \cX=(Q^o\backslash T) \times \PP^1$
the second factor can be equipped with
a homogeneous coordinate system $(\zeta_1: \zeta_2)$ such that the image, say, $\cD$
of $D^o=D\cap X^o$ in $X^o$ is defined by  the equation $\zeta_1=0$. We treat
$$
v=\zeta_1/\zeta_2
$$
as a coordinate in the neighborhood $\cX \backslash \{ \zeta_2 =0 \}$
of $\cD$ in $\cX$.

We fix a function $f \in \kk [Q]$ such that its pullback on $X$ belongs to $\ker \de
\backslash\ker \de_0$.
This pullback induces rational functions on $X^o$ and on $\cX$ denoted by the same symbol $f$.
By Lemma \ref{5.2}(1)
$f$ has poles along $D_0\cap X^o$.

By our choice of $B$ in Corollary \ref{5.6} $T^o$ is a submanifold of $Q^o$. Thus
locally the ideal of $T^o$ is generated by some function, say,  $u$ on $Q^o$. On $Q^o$ the function $f$
is of form $a/u^s$. Here $s\ge 1$ is the pole order  of $f$ along $T^o$, so $a$ is a rational function on $Q^o$,
which is nonzero in the general point of $T^o$.

Later on we will replace $f$ by a sufficiently large power $f^k$. By this we can achieve that the
pole order $s$ is arbitrary large.

Recall that  $U_f$ stands for the replica of $U$ associated with the locally nilpotent vector field $f\de$.
We note that  $U_f$ is well defined on the set
$$
X^o\backslash C \cong (Q^o\backslash T^o)\times \PP^1,
$$
cf.\ Corollary \ref{5.6}. Its element at moment $\tau \in \kk$ will be denoted by $h_{f, \tau}$.
Considered as an automorphism of $(Q^o\backslash T)\times \PP^1$ it
preserves the first factor but not the second one. The action of $h_{f, \tau}$ on $v$
is described
by the following Lemma.
\enota

\blem\label{5.8} There exist a regular function $d = d(f)$ on $Q^o$, which does not vanish at
 general points of $T$, and an
integer $l$ such that the automorphism of $(Q^o\backslash T)\times \PP^1$ defined by
$h_{f, \tau}$
is given in the coordinates
$(q,v)$ by
the formula
$$
h_{f,\tau}:\, (q,v)\mapsto \left( q\,\,,\, {\frac{u(q)^{m}v}{u(q)^{m} +\tau d(q) v}}\right)\,,
$$
where $m=s-l$. In particular $\cD \cap \cC= \{ u=v=0 \}$
is the set of indeterminacy points of $h_{f,\tau}$.
\elem

\bproof In homogeneous coordinates $(\zeta_1 : \zeta_2)$ the action of $U=\exp(\kk\de)$ on
$(Q^o\backslash T)\times\PP^1$
is of form $(\zeta_1 : \zeta_2) \to (\zeta_1, \zeta_2 + \tau c \zeta_1)$ where $c$
is non-vanishing function
on $Q^o\backslash T$.
That is, $c=c_0u^l$ where $c_0$ is a non-vanishing function on $Q^o$ and $l\in \ZZ$.
Hence $h_{f,\tau }$
is of form
$(\zeta_1 : \zeta_2) \mapsto (\zeta_1: \zeta_2 + {\frac{\tau d}{u^{s-l}}} \zeta_1)$,
where $d$ does not
vanish at general points  of $T^o$. Note that $m>0$ since $f\de$ has a pole along $D_0$. Passing to the
affine coordinate $v=\zeta_1/\zeta_2$ this  yields
the desired conclusion.
\eproof

Letting $s$ be the pole order of $f$ along $T$ we consider the set
\be\label{Pf}
P_f=\{q\in T:\mbox{ locally } f=a/u^s\mbox{ with } a(q)=0\mbox{ or }a\not\in\cO_{\bQ,q}\}\,,
\ee
where $u$ is as before (i.e.\ $u=0$ is a local equation of $T$ near $q$) and $a$ is a rational function.
This set is a proper closed subset of $T$.
The next proposition is the main result of this section.

\bprop\la{5.9}
Given  $m$ and a function $f\in \kk[Q]\cap \ker\de\backslash \ker \de_0$ there exists a positive integer $k_0$ such
that any transformation
$$
h\in U_{f^k}, \quad h\ne \id, \; k\ge k_0,
$$
is an $m$-contractions of $C$ along $D$ over the points of $Q^o\backslash P_f$.
\eprop

\bproof
Let $s,l$ be as in Notation \ref{5.7} and Lemma \ref{5.8}.
If we chose $k_0$ in such a way that $m'=k_0s-l\ge m$ then by Lemma \ref{x.11b} the map $h=h_{f^k,\tau}$
is indeed an $m$-contraction for any $\tau\ne 0$.
\eproof

Let now $Y\subseteq X$ be a closed subset. Consider the partial boundary
$$
\p_0Y=\bY\cap D_0\,.
$$
For $U\in \cU(X)$ we let $U^*=U\backslash \{\id\}$.
With this notation the following result holds.

\bprop\label{prop-one}
Let the notation and conventions be as in Notation \ref{5.1} and assume that (a), (b) in Lemma \ref{5.3} are satisfied.
Let $(Y_{\alpha,\beta})_{(\alpha,\beta)\in A\times B}$ be a flat family of proper closed subsets of $X$. Suppose
that there is a flat family $(E_\alpha)_{\alpha\in A}$ of proper, closed subset of $D$ such that
$$
\p Y_{\alpha,\beta}\cap D\subseteq E_\alpha \quad \mbox{for all}
\quad (\alpha,\beta)\in A\times B.
$$
Given an invariant function $f\in\ker\de\backslash\ker \de_0$, there is a dense open subset $A^o$ of $A$ and
a flat family $(E'_\alpha)_{\alpha\in A^o}$ of proper closed subset of $D_0$ satisfying
$$
\p_0 h.Y_{(\alpha,\beta)}\subseteq E'_\alpha
\quad
\forall\, (\alpha,\beta)\in A^o\times B,\,\,\forall\, h\in U^*_{f^k},\,\,\forall \, k\ge k_0\,.
$$
\eprop

\bproof
According to Proposition \ref{x.4} and Remark \ref{x.7} the closure $\bY_{\alpha\beta}$ of $Y_{\alpha\beta}$
in $\bX$ is at most $m$-tangent to $D$ along $C$ for $m\gg 0$ and for all $(\alpha,\beta)\in A\times B$ simultaneously.
Let $N_{\alpha\beta}=C\cap \ol{D\cap Y_{\alpha\beta}\backslash C}$ denote the defect set. By Proposition \ref{5.9}
for $k\gg 0$ any map $h\in U_{f^k}^*$ is an $m$-contraction of $C$ along $D$ over the points of $Q^o\backslash P_f$.
Applying Proposition \ref{x.13} the image $h.Y_{\alpha\beta}$ satisfies
\be\la{mmm}
\ol{ h.Y_{\alpha\beta}}\cap (D^o\cup C^o)\subseteq D\cup  h(N_{\alpha\beta})\cup  \brho^{-1}(P_f)
\ee
where $h(N_{\alpha\beta})$ stands for the total transform of $N_{\alpha\beta}$ under $h$.
By our assumption the defect set $N_{\alpha\beta}$ is contained in $N_\alpha =C\cap\ol{E_\alpha\backslash C}$.
Since our birational transformation $h$ is compatible with the fibration $\brho$, the total image $h(N_{\alpha\beta})$
is contained in $\brho^{-1}(\brho(N_\alpha))$. Taking in \eqref{mmm} the intersection with $D_0$ gives
$$
\p_0 (h.Y_{\alpha\beta})\subseteq E'_\alpha=(D\cup \brho^{-1}(B\cup \brho(N_\alpha)\cup P_f)) \cap D_0,
$$
where $B=\bQ\backslash \bQ^o$ is as in Corollary \ref{5.6}.
Using the theorem on generic flatness it is easily seen that  over an open dense subset $A^o$ of $A$ the sets
$E'_\alpha$ form a flat family of closed subsets of $D_0$.  This yields the assertion.
\eproof

\section{Proof of the main theorem}\label{sec4}

\subsection{Algebraic families of automorphisms}\label{transitivity results}

Following Ramanujam \cite{Ram} let us introduce the following notion.

\bdefi\label{algebraic family}  Given irreducible algebraic varieties $X$ and $A$ and
a map $\varphi:A\to\Aut(X)$ we say that $(A,\phi)$
is an {\em algebraic family of automorphisms on $X$} if the induced map
$A\times X\to X$, $(\alpha,x)\mapsto \varphi(\alpha).x$, is a morphism.
\edefi

By abuse of notation, we do not distinguish in the sequel  $A$ and its image $\varphi(A)$,  and we identify
$\alpha\in A$ with  its image $\varphi(\alpha)$ in $\Aut(X)$. As in the case of group action, given a point $x\in X$
the set $A.x$ will be called the {\em $A$-orbit} of $x$,
and the set $A_x=\{\alpha\in A\,|\,\alpha(x)=x\}$ the {\em stabilizer} of $x$ in $A$. The stabilizer admits
a natural linear representation
$d_x: A_x\to\GL(T_xX)$, $\alpha\mapsto d\alpha|T_xX$, called the tangent representation.

The following result allows to work with finite dimensional algebraic
families instead of dealing with infinite dimensional groups of automorphisms.

\blem\label{uuuuyo}
Let $X$ be a smooth quasi-affine variety and $G=G_{\cN}$ a group of automorphisms generated by a saturated
set of locally nilpotent derivations so that $G$ acts transitively on $X$. Then there exists
an algebraic family of automorphisms
$A\subseteq G$
such that  for any $x\in X$ we have

(a) $A.x=X$ and

(b) $d_x(A_x)=\SL(T_xX)$.
\elem

\bproof
According to Proposition 1.5 in \cite{AFKKZ} there exist one-parameter unipotent
subgroups $H_1,\ldots,H_s$ of $G$ such that with $H=H_1\cdot\ldots\cdot H_s\subseteq G$ we have
$H.x=G.x$ for any $x\in X$.
In particular, (a) holds with the algebraic family $A=H$.

By Theorem 4.2 \cite{AFKKZ} and its proof, for a fixed point $x\in
X$ the group $\SL(T_xX)$ is equal to the image in
$d_x(H')\subseteq \GL(T_xX)$ for an algebraic family
$H'=H'_1\cdot\ldots\cdot H'_r$, where $H_1',\ldots ,H_r'$ are
suitable one-parameter subgroups of $G_{\cN,x}$. Taking the
product $A=HH'H^{-1}$, where $H$ is as in (a) and
$H^{-1}=H_s\cdot\ldots \cdot H_1$,  we thus achieve that both (a)
and (b) are satisfied at every point  $x\in X$. \eproof

\bnota\label{44a} (a) As before we let $X$ be a smooth
quasi-affine variety and $G=G_{\cN}$ a group of automorphisms
generated by a saturated set of locally nilpotent derivations as
in Notation \ref{5.1}(a). We suppose that $G$ acts transitively on
$X$. According to Theorem \ref{1.14} there are derivations $\delta_0,\delta_1\in \cN$ such that the group
$$
H=\lla \de_0,\de_1\rra\subseteq G
$$
generated by $\de_0,\de_1$ and their replicas  acts with an open orbit on $X$.\footnote{In
contrast to Notation \ref{5.1}(a)
in this section the role of $\de_0$ and $\de_1$ will be symmetric
so that it is convenient to replace the former
$\de$ by $\de_1$.}
These locally nilpotent vector fields generate  one-parameter unipotent subgroups
$U^0, U^1\in\cU(G)$.
Any function $f\in \ker\de_0\backslash \ker \de_1$ yields a replica $U^0_f$, and similarly
$g\in \ker\de_1\backslash \ker \de_0$ yields a replica $U^1_g$.

(b) To any sequence of invariant functions
\be\label{seq}
\cF=\{f_1,\ldots,f_s, g_1, \ldots,g_s\},\,\,\,\,\mbox{where}
\,\,\,\, f_i\in\ker\de_1\backslash\ker
\de_0\,\,\,\,\mbox{and}\,\,\,\, g_i\in\ker\de_0\backslash\ker
\de_1\,,
\ee
we associate an algebraic family
of automorphisms $\AA^{2s}\to\Aut(X)$ defined by the product
\be\label{00121}
U^\cF=U^1_{f_s}\cdot
U^0_{g_s}\cdot\ldots\cdot U^1_{f_1}\cdot
U^0_{g_1}\subseteq H\,.
\ee
More generally, given a tuple
$\kappa=(k_i,l_i)_{i=1,\ldots,s}\in\N^{2s}$ the product
\be\label{001210}
U_\kappa=U_\kappa^\cF=
U^1_{f_s^{k_s}}\cdot
U^0_{g_s^{l_s}}\cdot\ldots\cdot U^1_{f_1^{k_1}}\cdot
U^0_{g_1^{l_1}}
\subseteq H\,
\ee
is as well an algebraic family of automorphisms.
\enota

\bcor\label{cor-0012} There is a finite collection of invariant functions $\cF$ as in (\ref{seq})
such that for any sequence
$\kappa=(k_i,l_i)_{i=1,\ldots,s}\in\N^{2s}$ the algebraic family of automorphisms
$U_\kappa$ as in (\ref{001210})
has a dense open orbit in $X$.
This orbit $O(U_\kappa)$ coincides with $O(H)$ and so
does not depend on the choice of $\kappa\in\N^{2s}$.\ecor

\bproof According to Proposition 1.5 in \cite{AFKKZ} there is a
sequence $\cF$ as in (\ref{seq}) such that
$$
H.x= U^\cF.x\quad\forall x\in X
\,.
$$
In
particular, for $x\in O(H)$ the orbit $ U^\cF.x=O(H)$ is open in
$X$. It is easily seen that for any $\kappa\in\N^{2s}$ we have
$O(U_\kappa)=O( U^\cF)=O(H)$. Indeed,  $O(H)$ consists of all the $U^\cF$-flexible points in $X$.
Now the assertions follow. \eproof

\subsection{Proof of the main theorem}

\bnota\label{4.5}
We keep the notation and assumptions from \ref{44a}(a).

(a) Let $\rho_0:X\to Q_0$ and $\rho_1:X\to Q_1$ be partial quotients with respect
to the unipotent subgroups $U^0$ and $U^1$, respectively.
Let us choose open embeddings $X\hto \bX$, $Q_0\hto\bQ_0$, and $Q_1\hto\bQ_1$
into normal projective varieties, see
Notation \ref{5.1}. We can assume that the following conditions
are satisfied.
\bnum[(i)]
\item  $\rho_0$ and $\rho_1$ extend to morphisms
$\brho_0:\bX\to\bQ_0$ and $\brho_1:\bX\to \bQ_1$.
Let $D_0$ and  $D_1$ as in \ref{5.1} be the unique horizontal divisors that map birationally onto
$\bQ_0$ and $\bQ_1$, respectively.

\item $\bX$, $D_0$ and $D_1$ are smooth, see Lemma \ref{5.3}(b).

\item $T_0=\brho(D_0)$ and $T_1=\brho(D_1$ are divisors in $\bQ_0$ and $\bQ_1$, respectively;
see Lemma \ref{5.3}(a).
\enum
(b) Given a closed subscheme $Y\subseteq X$ of codimension $\ge 2$ we call
$$
\p_0Y=\bY\cap D_0\and \p_1Y=\bY\cap D_1
$$
the {\em partial boundaries.} Furthermore $O_Y$ will denote the open orbit of $G_{\cN,Y}$ in $X\backslash Y$.
\enota

\bsit\label{cvnt}
In the course of the proof of the main Theorem we move the given pair $(Y,x)$ to
another one $(Y_\alpha, x_\alpha)$ by means of an automorphism $\alpha \in G_\cN$, where
$Y_\alpha=\alpha. Y$ and $x_\alpha=\alpha.x$. In this way we can adopt
the position of our pair with respect to the $\PP^1$-fibration
$\brho_0:\bX\to\bQ_0$
so that the conditions  (i)-(iii) below hold.

\bnum\la{condis}
\item[(i)] $U^0.x_\alpha\cap O_{Y_\alpha}\neq\emptyset\,;$
\item[(ii)] $U^0.x_\alpha\cap Y_\alpha=\emptyset\,;$
\item[(iii)] $\p_0(U^0.x_\alpha)\notin\p_0 (Y_\alpha)$.\enum
\esit

The following lemma allows to deduce Theorem \ref{mthm1}
provided that (i)-(iii) hold for any $x\in X\backslash Y$ with some $\alpha\in G$ depending on $x$.

\blem\label{777} If for a point $x\in X\backslash Y$ and for some $\alpha\in G$ conditions (i)-(iii) in \ref{condis}
are fulfilled then $x\in O_{Y}$.
If these conditions are fulfilled for any $x\in X\backslash Y$ with some $\alpha\in G$ depending on $x$,
then the conclusion of Theorem  \ref{mthm1} holds.
\elem

\bproof
Since $O_{Y\alpha}=\alpha.O_{Y}$ we have
$$
x\in O_{Y}\Longleftrightarrow x_\alpha\in O_{Y_\alpha}\,.$$
Replacing $(Y,x)$ by $(Y_\alpha,x_\alpha)$ we will assume that (i)-(iii) hold for the pair
$(Y,x)$ and $\alpha=\id$. We need to show that then $x\in O_Y$. Conditions  (ii) and (iii) yield that
$$
\rho_0(x)\in\rho_0(O_{Y})\backslash\overline{\rho_0(Y)}\,.
$$
Therefore there exists a regular function
$h\in\cO(Q_0)$ such that $h(\rho_0(x))=1$ and $h$ vanishes
on $\rho_0 (Y)$. Replacing $h$ by  a suitable power of $h$
we may suppose that the $\de_0$-invariant function $f=h\circ\rho_0$ on $X$ vanishes on $Y$.
Thus the replica $U^0_f=\exp(\kk f\de_0)$ of $U^0$ fixes $Y$ pointwise i.e.
$U^0_f\in \cU(G_{\cN,Y})$.
By (i) one can find $u\in U^0_f$ such that $u.x\in
O_{Y}$. Hence also $x\in O_{Y}$, as stated.
\eproof

Thus to prove Theorem \ref{mthm1} it is enough to show that (i)-(iii) hold for every point
$x\in X\backslash Y$ with a suitable $\alpha\in G$ depending on $x$.

\blem\label{constr} Given a point $x\in X\backslash Y$ and
an algebraic family  of automorphisms  $\varphi:A\to\Aut(X)$ the following hold.

(a) The set of all $\alpha\in A$ satisfying (i)  is open in $A$.

(b) The set of all $\alpha\in A$ satisfying  (ii)  is constructible in $A$.
\elem

\bproof
(a) The subset $B\subseteq A$ where (i) does not hold is the set of $\alpha\in A$ satisfying
$$
U^0.x_\alpha\subseteq Y_\alpha
\quad\mbox{or, equivalently,}\quad \alpha^{-1}U^0\alpha. x\subseteq Y.
$$
Thus $B=\bigcap_{u\in U^0} B_u$, where
$B_u=\{\alpha\in A: \alpha^{-1}u\alpha.x\in Y\}$
 is the preimage of $Y$ under the morphism $A\to X$, $\alpha\mapsto
\alpha^{-1}u\alpha.x$. Hence $B$ is closed in $A$. This proves (a).

(b) Similarly, the subset $C\subseteq A$ where (ii) does not hold is the set of $\alpha\in A$ with
$\alpha^{-1}U^0\alpha\cap Y\ne \emptyset.$
Consider the set
$$
C'=\{(\alpha, u)\in A\times U^0: \alpha^{-1}u\alpha.x\in Y\}\,.
$$
This set is closed in $A\times U^0$ since it is the preimage of $Y$
under the morphism $A\times U^0\to X$,  $(\alpha,u)\mapsto
\alpha^{-1}u\alpha.x$. Since $C$ is the image of $C'$ under the projection to $A$, (b) follows.
\eproof

The next proposition allows to verify conditions (i) and (ii).

\bprop\label{555} Let as before  $x\in X\backslash Y$.
\bnum[(a)]
\item If $A$ is an algebraic family of automorphisms of $X$ with $d_x(A_x)\supseteq\SL(T_xX)$,
then the set of all $\alpha\in
A$ satisfying (i) is a dense open subset of $A$.

\item
There exists an algebraic family $A^*\subseteq G_x$
transitive in $X^*=X\backslash\{x\}$ such that for any subgroup $U^0\in \cU(X)$ condition
(ii) holds for a general $\alpha\in A^*$.

\item
Given an algebraic family $B\subseteq\Aut(X)$
we let  $\tilde A=B\cdot A^*\subseteq \Aut(X)$,
where $A^*\subseteq G_x$
is as in (b).
Then  (ii) holds for a general $\tilde\alpha\in \tilde A$.
\enum\eprop

\bproof (a) By Lemma \ref{constr} it suffices to find $\alpha\in A$ satisfying (i),
or, equivalently, such that $\alpha^{-1}U^0\alpha.x\cap O_{Y}\neq\emptyset$.
By our  assumptions in (a) for any nonzero vector $v\in T_xX$ there is an element $\alpha\in A_x$ such that
$v$ is tangent to the orbit through $x$ of the one-parameter group $\alpha^{-1}U^0\alpha\subseteq\Aut(X)$.
These orbits form
an algebraic family of smooth rational curves in $X$ through the point $x$ that dominates $X$ and so
meets the open orbit $O_{Y}$, as required.

(b)  By the Transversality Theorem  \cite[1.16]{AFKKZ} there exists an algebraic family $A^*\subseteq G_x$
transitive in $X^*$
such that for any two subvarieties $Y,Z\subseteq X$ there is a  dense open subset
$A_0\subseteq A^*$ with the property that
for any $\alpha\in A_0$ the varieties $\alpha.Y$ and $Z$ are transversal. Applying this to $Z=U^0.x$
the varieties $U^0.x$ and $\alpha.Y$ are disjoint,  because under our assumptions
$$\dim U^0.x +\dim Y< \dim X\,.$$
Since $x_\alpha=x$, (b) follows.

To deduce (c) we note that the set, say $C$ of points $\tilde \alpha\in \tilde A$, where (ii) fails is the set of
$\tilde\alpha=(\beta,\alpha)$ with $\alpha^{-1}\beta^{-1}U^0\beta\alpha.x\cap Y\ne\emptyset$. Consider
similarly as in the proof of Lemma \ref{constr}(b) the closed subset of $B\times A^*\times U^0$
$$
C'=\{(\beta,\alpha,u)\in B\times A^*\times U^0: \alpha^{-1}\beta^{-1}u\beta\alpha.x
\in  Y\}\,,
$$
where $A^*$ satisfies the conclusion of (b).
According to (b) for any $\beta\in B$ the  set
$$
C'_\beta=C'\cap (\{\beta\}\times A^*\times U^0)
$$
maps under the projection to $A^*$ to a nowhere dense subset.
Hence also the image $C$ of $C'$ under the projection to $\tilde A=B\times A^*$ will be nowhere dense.
Thus its complement contains an open dense subset proving (c).
\eproof

\bnota\label{nota4.10}
Given a one-parameter group $U\in \cU(X)$ we let as before $U^*=U\backslash\{\id\}$.
Given a collection $\cF$ of invariant functions
$$
f_1,\ldots,f_s\in\ker\de_1\backslash\ker \de_0\quad\mbox{and}\quad g_1,
\ldots,g_s\in\ker\de_0\backslash\ker \de_1
$$
and $U_\kappa=
U^1_{f_s^{k_s}}\cdot
U^0_{g_s^{l_s}}\cdot\ldots\cdot U^1_{f_1^{k_1}}\cdot
U^0_{g_1^{l_1}}$ as in \eqref{00121},  we let
$$
U_\kappa^*=
U^{1*}_{f_s^{k_s}}\cdot
U^{0*}_{g_s^{l_s}}\cdot\ldots\cdot U^{1*}_{f_1^{k_1}}\cdot
U^{0*}_{g_1^{l_1}}\,.
$$
\enota
Using Proposition \ref{prop-one} we can deduce the following result.

\bprop\label{prop-two}
Let $(Y_\alpha)_{\alpha\in A }$ be a flat family of proper closed subsets of $X$.
Assume that the partial boundaries
$\p_i Y_\alpha$ (see Notation \ref{4.5}) are contained in $E_{\alpha,i}$,
where the $(E_{\alpha,i})_{\alpha\in A}$, $i=0,1$,
form flat families of proper closed subsets of $D_i$.
Then one can find an open dense subset $A^o$ of $A$, flat families of proper, closed subsets
$(E^o_{\alpha, i})_{\alpha\in A^o}$ of $ D_i$ ($i=0,1$), and a sequence
$\kappa=(k_1,l_1,\ldots,k_s,l_s)\in\N^{2s}$ such that for any element $h\in U_\kappa^*$ we have
$$
\p_i (h.Y_\alpha)\subseteq E^o_{\alpha ,i}\, ,\qquad i=0,1\,,\; \forall\,\alpha\in A^o\,.
$$
\eprop

\bproof
The proof proceeds by induction on $s$. For $s=0$
the assertion clearly holds with $A^o=A$ and
$E_{\alpha, i}=\p_i Y_\alpha$, $i=0,1$. Assume that it holds at step $s-1$, i.e. we can find
$\kappa'=(k_j,\, l_j)_{j=1,\ldots,s-1}\in\N^{2s-2}$, a dense open subset $A'\subseteq A$
and flat families of proper closed subsets  $(E_{\alpha,i})_{\alpha\in A'}$ of $ D_i$ such that for $\alpha\in A'$
$$
\p_i(h.Y_\alpha)\subseteq E_{\alpha,i}\,,\quad i=0,1,\quad \forall\,
h\in U_{\kappa'}^*\,.
$$
The varieties $(h.Y_\alpha)_{(h,\alpha)\in U_{\kappa'}^*\times A'}$ form a flat  algebraic family.
By Proposition \ref{prop-one}
one can find an open dense subset $A''\subseteq A'$ and flat families $(E'_{\alpha,i})_{\alpha\in A''}$, $i=0,1$,
of proper closed subsets of $D_i$ such that
$$
\p_i(h'h.Y_\alpha)\subseteq E'_{\alpha,i}\,\,\, (i=0,1)\quad
\forall \,\, l_s\gg 0,\,\, \forall\alpha\in A'',\,\,
\,\,\forall \,  (h',h)\in U^{0*}_{g_s^{l_s}}\times U_{\kappa'}^*\,.
$$
Fixing a sufficiently large $l_s$ and applying the same argument again
one can find an open dense subset $A^o\subseteq A''$ and flat families $(E^o_{\alpha,i})_{\alpha\in A^o}$, $i=0,1$,
of proper closed subsets  of $D_i$ such that
$$
\p_i(h''h'h.Y)\subseteq E^o_{\alpha,i}\,\,\, (i=0,1)\quad
\forall k_1\gg 0, \,\, \forall\alpha\in A^o\,,\,
\forall  \,  (h'', h',h)\in  U^{1*}_{f_s^{k_s}}\times U^{0*}_{g_s^{l_s}}\times U_{\kappa'}^*\,.
$$
This concludes the induction.
\eproof

Using Proposition \ref{prop-two} and
Corollary \ref{cor-0012} we can now
deduce Theorem \ref{mthm1}.

\bproof[Proof of Theorem \ref{mthm1}.]
Let $x\in X\backslash Y$ be a fixed point.
We show that for a suitable choice of an algebraic family $A$ of automorphisms conditions (i)-(iii) are satisfied for
the pair $(Y_\alpha, x_\alpha)$, if $\alpha\in A$ is generic. Then our theorem follows by applying Lemma \ref{cvnt}.

{\em Step 1.} Consider an algebraic family $A\subseteq G$ satisfying conditions (a)
and (b) of Lemma \ref{uuuuyo}.
Applying Proposition \ref{555}(a) condition (i) holds
when $\alpha$ varies in a dense open subset of $A$.
Replacing the original pair $(Y,x)$ by a suitable new one
$(Y_\alpha,x_\alpha)=(\alpha.Y,\alpha.x)$ we may suppose that $(Y,x)$ satisfies (i).

{\em Step 2.}
In the following we construct an algebraic family $B$ of automorphisms such that for a generic choice of $\beta\in B$
the translates $(Y_\beta, x_\beta)$ satisfy (ii), (iii). Since by Proposition \ref{555}(a) condition (i) is open then the pair
$(Y_\beta, x_\beta)$ also satisfies (i).

Let $A^*$ be a family of automorphisms as in Proposition \ref{555}(b). The translates $Y_\alpha=\alpha. Y$,
$\alpha\in A^*$, form a flat family of proper closed subsets of $X$. Using the theorem of generic flatness it is
easily seen that over an open dense subset $A'\subseteq A^*$ also the partial boundaries $E_{\alpha,i}=\p_i Y_\alpha$,
$\alpha\in A'$, form flat families of proper closed subsets of $D_i$, $i=0,1$.
Let now $\cF$, $U_\kappa$,
and $U_\kappa^*$ be as in Notation \ref{nota4.10}.
By Proposition \ref{prop-two} we can find $\kappa=(k_1,l_1,\ldots,k_s,l_s)\in\N^{2s}$,
a dense open subset $A^o\subseteq A'$,
and families $(E^o_{\alpha,i})_{\alpha\in A^o}$, $i=0,1$, of proper closed subsets of $D_i$ such that
$\p_i(h.Y_\alpha)\subseteq E^o_{\alpha, i}$ for $i=0,1$, $\alpha\in A^o$ and all $h\in U_\kappa^*$.

We claim that for a generic choice of $(h,\alpha)\in B=U^*_\kappa\times A^*$ conditions (ii) and (iii) are
satisfied for $h.Y_\alpha$. To check (ii) we note that $h.Y_\alpha=h\alpha.Y$. Thus applying Proposition \ref{555}(c)
to the family $B=U^*_\kappa\times A^*$ condition (ii) is indeed satisfied for a generic choice of $(h,\alpha)$.

It remains to show that (iii) is satisfied for a generic choice of $(h,\alpha)$.
Condition (iii) is equivalent to $\brho_0(h.x_\alpha)\not \in \p_0(h.Y_\alpha)$. By construction
$\p_0(h.Y_\alpha)\subseteq E_{\alpha,0}\subseteq D_0$ for any $h\in U_\kappa^*$,
while for a fixed $\alpha\in A^o$ the points
$h.x_\alpha$, $h\in U_\kappa^*$, fill in a dense subset of $X$,
and so their images $\brho_0(h.x)$ fill in a dense subset of $Q_0\subseteq \bar Q_0\backslash \brho_0(D_0)$.
Thus (iii) holds for a generic choice of $(h,\alpha)\in U^*_\kappa\times A^o$.
This concludes the proof of Theorem \ref{mthm1}.
\eproof

\end{document}